\documentclass[11pt,letterpaper]{amsart}
\usepackage{arxiv_math}
\usepackage{import}
\usepackage{adjustbox}

\title{For exotic surfaces with boundary, one stabilization is not enough}
\author{Gary Guth} \address{University of Oregon, Eugene, OR, 97403}
\email{gguth@uoregon.edu}

\begin{document}
\begin{abstract}
    A result of Baykur-Sunukjian \cite{baykur_sunukjian_knotted_surf_stab} states that homologous surfaces in a 4-manifold become isotopic after a finite number of internal stabilizations, i.e. attaching tubes to the surfaces. A natural question is how many stabilizations are needed before the surfaces become isotopic. In particular, given an exotic pair of surfaces, is a single stabilization always enough to make the pair smoothly isotopic? We answer this question by studying how the stabilization distance between surfaces with boundary changes with respect to satellite operations. Using a range of Floer theoretic techniques, we show that there are exotic disks in the four-ball which have arbitrarily large stabilization distance, giving the first examples of exotic behavior in the four-ball for which ``one is not enough''.
\end{abstract}

\maketitle 

\tableofcontents

\section{Introduction}\label{section: intro}

The study of exotic phenomena in dimension four has been an active area of mathematical research for decades. A pair of smooth, closed 4-manifolds $X_1$ and $X_2$ are called \emph{exotic} if the two are homeomorphic but not diffeomorphic. A classic result of Wall \cite{Wall_simply_connected} shows that exotic behavior is ``unstable'', in the following sense: if $X_1$ and $X_2$ are an exotic pair of simply connected 4-manifolds, $X_1\#^n S^2\times S^2$ are diffeomorphic $X_2\#^n S^2\times S^2$ for sufficiently large $n$. This fact follows from a simple observation: two smooth, simply connected, oriented 4-manifolds which are homeomorphic are cobordant through a cobordism with only 2- and 3-handles. An intermediary slice of this cobordism can be described as \emph{either} $X_1\#^n S^2\times S^2$ or as $X_2\#^n S^2\times S^2$. However, as this cobordism is constructed abstractly, it is completely opaque how many $S^2\times S^2$ summands are needed before the exotic behavior dissolves. Work of Gompf \cite{Gompf1984StableDO} extends this result for oriented manifolds with arbitrary fundamental group (in the non-orientable case, it is necessary to allow connected sums with $S^2\tilde{\times} S^2$ as well).

An outstanding conjecture is that, in the oriented case, a single stabilization (taking a connected sum with a single $S^2\times S^2$) is always enough to eliminate any exotic behavior. To the author's knowledge, every known construction of exotic oriented, closed 4-manifolds yields pairs which become diffeomorphic after a single stabilization \cite{baykur_dissolving_knot_surgered,auckly_families_one_stab,akbulut_variations_FS_knot_surgery,baykur_sunukjian_round_handles_log_transforms}.  \\

A related question asks whether homologous surfaces with the same genus in given 4-manifold are always ``stably equivalent''. Two surfaces are called exotic if they are topologically isotopic, but not smoothly isotopic. In the surface case, there are two natural notions of stable equivalence. By work of Quinn \cite{quinn_iso_4mflds} and Perron \cite{perron_pseudo_isos, perron_pseudo_iso_and_isos} homologous surfaces with equal genus in a 4-manifold become isotopic after sufficiently many \emph{external stabilizations}, connected sums with $S^2 \times S^2$ away from the surfaces. Recent work of \cite{AKMRS_isotopy_surfaces_single_stab} gives a partial answer to the analogous ``one is enough'' conjecture in this context, by showing that under some algebraic-topological conditions, a single external stabilization is enough to make homologous surfaces with equal genus isotopic. Baykur-Sunukjian give another version of stability for surfaces, showing that any pair of homologous surfaces become isotopic after a sufficient number of \emph{internal stabilizations}, i.e. by attaching unknotted 1-handles to the surfaces. In addition, they prove for every known construction of pairs of exotic surfaces, a single internal stabilization is enough to dissolve the exotic behavior \cite{baykur_sunukjian_knotted_surf_stab}. \\

However, the positive evidence for the ``one is enough'' conjecture in the closed case is countered by negative evidence in the relative case. Miller-Powell \cite{MP_stab_dist} and Juh\'asz-Zemke \cite{juhasz_zemke_stabilization_bounds} produce disks with common boundary a knot in $S^3$ which remain distinct after many internal stabilizations (though, their examples are topologically distinct, as well as smoothly). Lin-Mukherjee \cite{LM_FamilyBaurFuruta_Exotic_Surf_Smale} produce an exotic pair of links in a punctured $K3$-surface which remain exotic after performing an external stabilization. In this paper, we produce the first examples of exotic surfaces in the four-ball for which one internal stabilization is not enough.

\begin{theorem}\label{Theorem 1: exotic disks with large stabilization distance}
    For any $p$, there exists a knot $J_p$ which bounds a pair of exotic disks $D_p$ and $D_p'$ which remains exotic after $p-1$ internal stabilizations.
\end{theorem}

\subsection{Stabilization distance}\label{subsection: intro, stab dist}

Following \cite{MP_stab_dist} and \cite{juhasz_zemke_stabilization_bounds}, we approach this question by studying the \emph{stabilization distance} between two disks. 

\begin{defn}\label{def: preliminary def of stabilization distance}
Let $\Sigma$ and $\Sigma'$ be a pair of surfaces in $B^4$ with equal genera with boundary $K \sub S^3$. The \emph{1-handle stabilization distance}, $d(\Sigma, \Sigma')$ is the minimal number of 1-handle stabilizations needed to make $\Sigma$ and $\Sigma'$ isotopic. For a more precise statement, see Definition \ref{def: 1-handle stabilization}.
\end{defn}

\begin{figure}[h!]\centering
    \includegraphics[scale=.65]{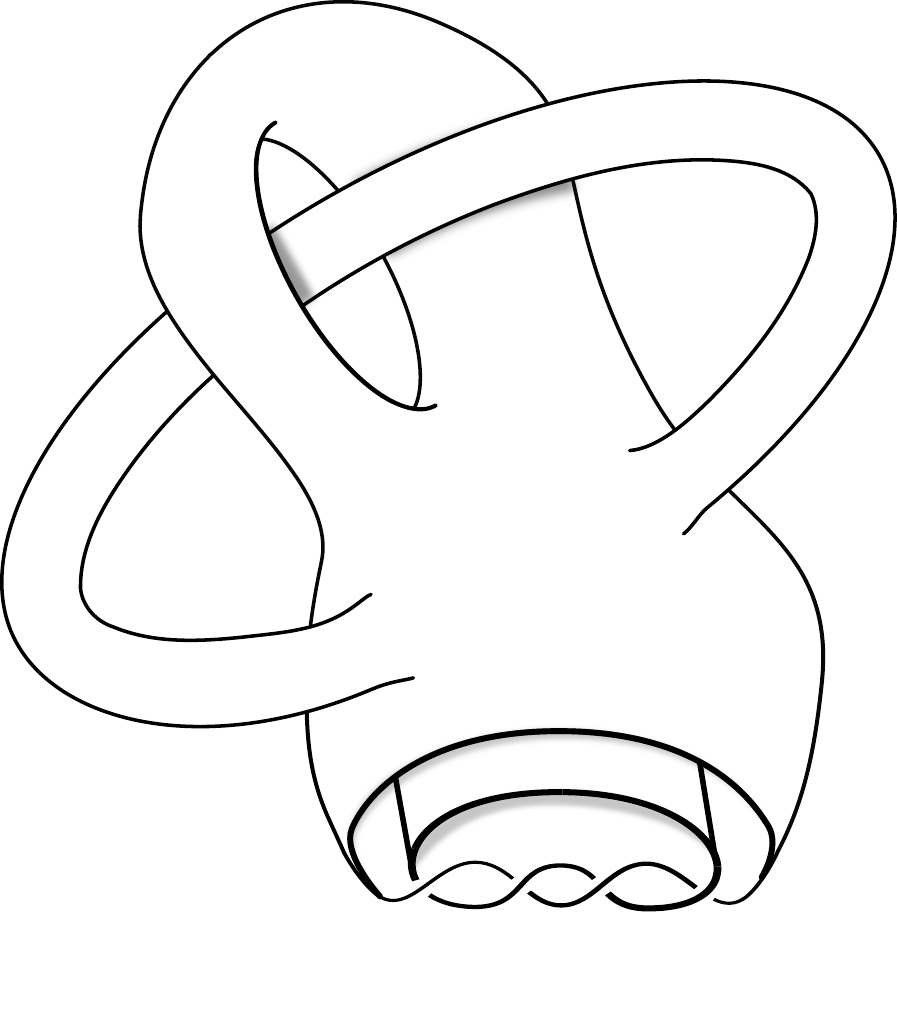}
    \caption{A 1-handle stabilization of a surface $\Sigma$.}
    \label{fig: 1-handle stabilization}
\end{figure}

Theorem 1 states $d(D_p, D_p') \ge p$, but we show somewhat more. We show that, in fact, $d(D_p, D_p') = p$. Additionally, we make use of a more general notion of stabilization distance, due to \cite{juhasz_zemke_stabilization_bounds}, and show our examples have large stabilization distance in that sense as well. For more details, see Section \ref{Section: Preliminaries}. \\

We distinguish our surfaces by comparing their induced maps on knot Floer homology. Stabilization has a simple effect on the induced map; attaching a tube simply corresponds to multiplication by $U$ (or $V$) \cite[Lemma 3.3]{juhász_zemke_newslicegenus},\cite{JMZ_TorsionOrder}. Juh\'asz-Zemke make use of this fact to define a suite of ``secondary'' Heegaard Floer invariants which provide lower bounds for the stabilization distance of two surfaces. 

Recall that the \emph{torsion order} $\Ord_U(K)$ of a knot $K$ is the smallest power of $U$ which annihilates the torsion submodule of $\HFK^-(K)$. Since a stabilization corresponds to multiplication by $U$, any two maps induced by disks with boundary $K$ become indistinguishable after multiplication by $U^{\Ord_U(K)}$. Therefore, these lower bounds cannot be used to show disks bounding knots with torsion order 1 have large stabilization distance. $\Ord_U(K)$ is bound above by the \emph{fusion number} of $K$, which is the minimal number of bands occurring in ribbon disk for $K$ \cite{JMZ_TorsionOrder}. Therefore, to have any hope in finding disks with large stabilization distance, it is necessary to work with knots with large fusion number.

\subsection{Cabled concordances} \label{subsection: intro, cabled concordances}

Recent work of Hom-Kang-Park \cite{HKP_ribbon_cables} and Hom-Lidman-Park \cite{HLP_unknotting_cables} studies how cabling is related to the torsion order and fusion number of a knot. If $K$ is ribbon with fusion number 1, then the knot Floer homology of the $(p, 1)$-cable of $K$ has torsion order $p$ \cite[Lemma 3.3]{HKP_ribbon_cables}. Cabling has a natural four-dimensional extension: given a concordance $C: K \ra K'$, there is a cabled concordance between the cables of $K$ and $K'$. In particular, given a ribbon knot $K$ with fusion number 1, $K_{(p, 1)}$ bounds a ``cabled'' ribbon disk, and has fusion number $p$. \\ 

\begin{figure}[h!]
    \includegraphics[scale=.4]{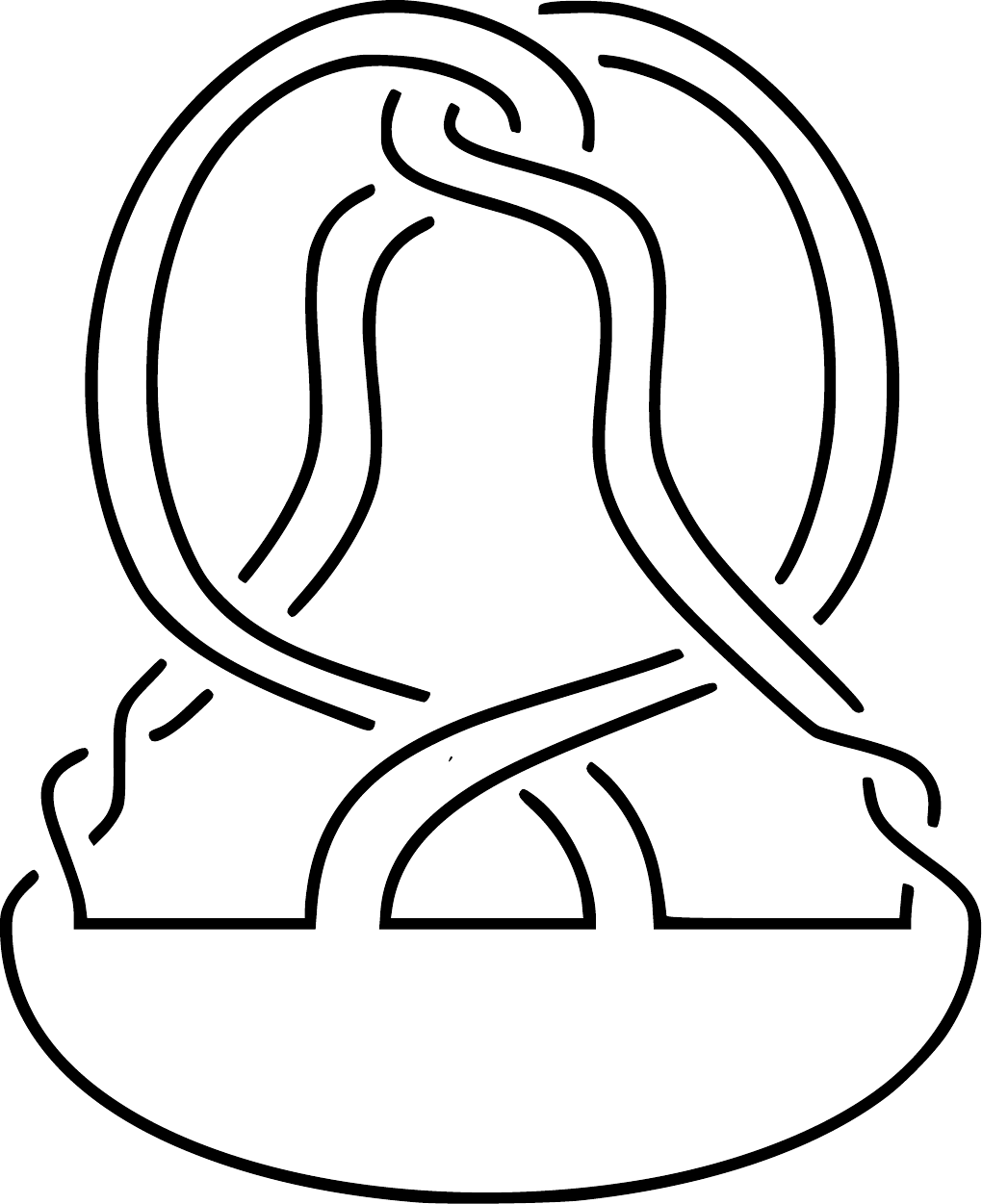}
    \caption{The knot $J$ which bounds an exotic pair of disks.}
    \label{fig: the knot J}
\end{figure}

The knot $J$ shown in Figure \ref{fig: the knot J}, bounds an exotic pair of disks $D$ and $D'$ by the work of Hayden \cite{Hayden_CorksCoversComplex}. Moreover, these disks are distinguished by their induced maps on knot Floer homology \cite{DaiMallickStoffregen_equivar_knots_floer}. But, $J$ has fusion number 1, so the two maps become equal after a single stabilization. In fact, we will show directly that these two disks are smoothly isotopic after a single stabilization. However, the cabled disks $D_p$ and $D_p'$ have fusion number $p$, and as we show, have stabilization distance $p$ as well.

In order to prove Theorem \ref{Theorem 1: exotic disks with large stabilization distance}, we show that the map on knot Floer homology induced by the cabled disk is determined by the map induced by the original disk. We show that, in general, the concordance induced by any satellite pattern can be computed in terms of the original concordance. 

\begin{theorem}\label{Theorem 2: satellite formula for concordance maps}
    Let $C: K \ra K'$ be a smooth concordance. Then, there exists a map $F: \CFDh(S^3-K) \ra \CFDh(S^3-K')$ induced by $C$, such that for any pattern knot $P$ in the solid torus, the following diagram commutes up to homotopy:

    \begin{center}
        \begin{tikzcd}
            \CFA^-(\cH_P) \boxtimes \CFDh(S^3-K) \ar[d, "\id \boxtimes F"]  \ar[r, "\simeq"]  
                & \CFK^-(K_{P}) \ar[d, "F_{C_P}"]  \\
            \CFA^-(\cH_P) \boxtimes \CFDh(S^3-K') \ar[r, "\simeq"]
                & \CFK^-(K_{P}'),
        \end{tikzcd}
    \end{center}
where $\cH_P$ doubly pointed, bordered Heegaard diagram for $P\sub S^1\times D^2$, $K_P$ and $K_P'$ are satellites of $K$ and $K'$, and $C_P$ is the concordance induced by $P$. The horizontal arrows are given by the pairing theorem \cite{LOT_bordered_HF}. 
\end{theorem}

As another application of Theorem \ref{Theorem 2: satellite formula for concordance maps}, we produce a collection of pairs of exotic cobordisms $W_n, W_n': L(n, 1) \ra S_n^3(J)$ by doing surgery on the concordances $C$ and $C'$. We distinguish these 4-manifolds by directly computing the maps $\wh{F}_{W_n}$ and $\wh{F}_{W_n'}$. 

\begin{theorem}\label{theorem: exotic cobordisms}
    The are exotic pairs of cobordisms $W_n, W_n': L(n, 1) \ra S_n^3(J)$ which are distinguished by their induced maps on $\HFh$.
\end{theorem}

\subsection{Acknowledgements} The author would like to thank his advisor, Robert Lipshitz, for his guidance and many helpful suggestions, and in particular, sharing the proof of Lemma \ref{lemma: closest point map}. Thanks, also, to Holt Bodish, Jesse Cohen, Maggie Miller, Ina Petkova, and Ian Zemke for helpful comments and conversations. Thanks to Arunima Ray for pointing out a mistake in the introduction. GG was partially supported by NSF Grant DMS-220214.

\section{Preliminaries}\label{Section: Preliminaries}

We begin by reviewing the two definitions of stabilization distance and illustrate that the two need not agree. We then review the construction of cabled concordances and show how this operation can be used to produce disks which are topologically isotopic. We conclude this section by giving an upper bound on the stabilization distance of $D_p$ and $D_p'$ by explicitly showing they become isotopic after $p$ stabilizations.

\subsection{Two notions of stabilization distance}\label{subsection: prelims, two defs of stabilization}

The most general notion of internal stabilization is due to \cite{juhasz_zemke_stabilization_bounds}. 

\begin{figure}[h!]\centering
    \includegraphics[scale=.4]{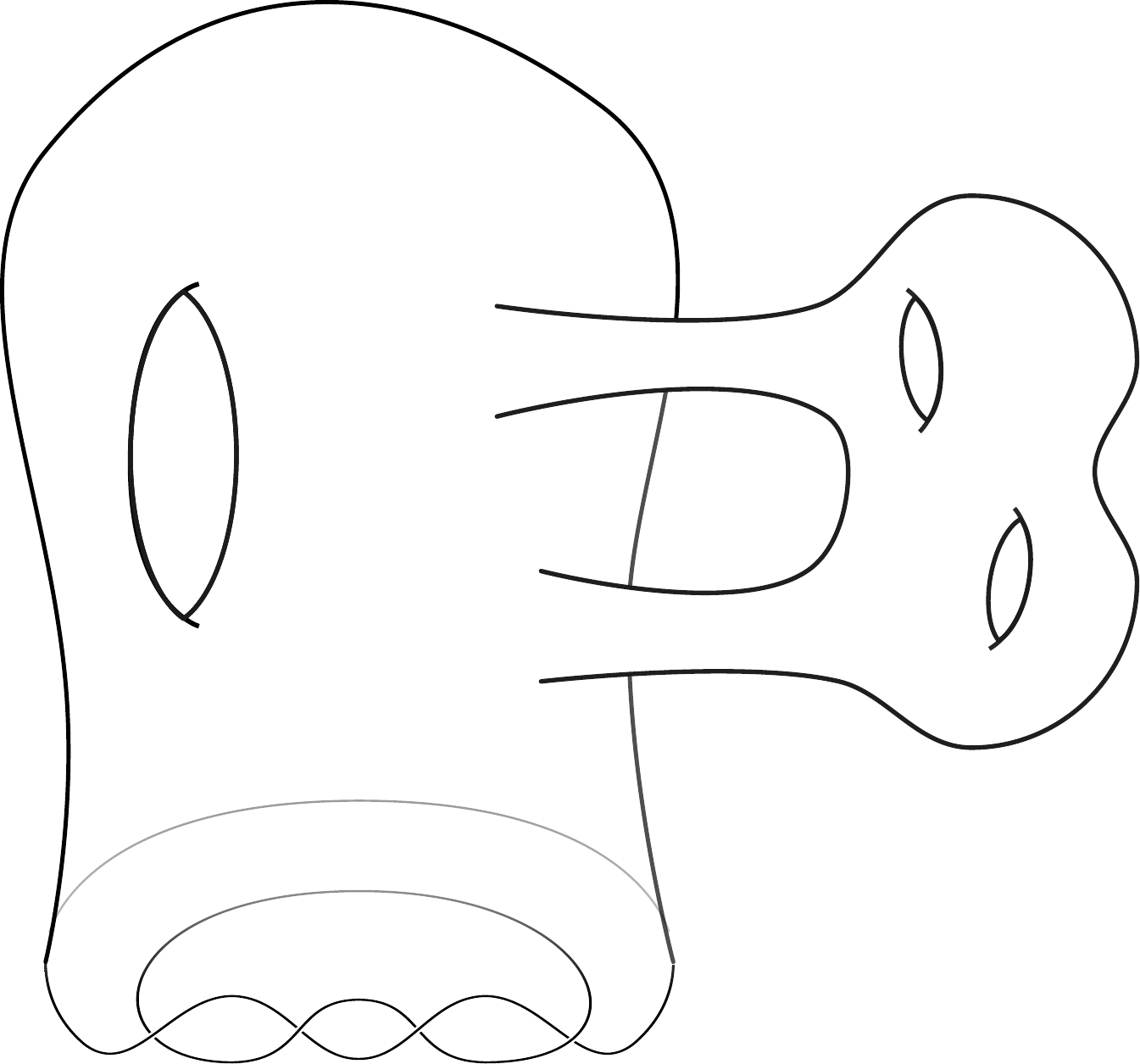}
    \caption{A $(g, n)$-stabilization along $(B^4,S_0)$. The case $(g, n)=(2, 2)$ is shown.}
    \label{fig: gn stab stabilization}
\end{figure}

\begin{defn}\label{def: (g, n) stabilization}
    Let $\Sigma$ be an oriented surface with boundary, smoothly embedded in a 4-manifold $W$. Let $B$ be a 4-ball in the interior of $W$ whose boundary intersects $\Sigma$ in an $n$-component unlink $L$. Moreover, suppose $\Sigma \cap B$ is a collection of disks $D_1, ..., D_n$ which can be isotoped into $\partial B$ relative to their boundaries. Let $S_0$ be a connected genus $g$ surface in $B$ with boundary $L$. The surface $\Sigma' = (\Sigma - B) \cup_L S_0$ is called the \emph{$(g, n)$-stabilization of $\Sigma$ along $(B, S_0)$}. We call $\Sigma$ the \emph{$(g, n)$-destabilization of $\Sigma'$ along $(B, S_0)$}. See Figure \ref{fig: gn stab stabilization}.
\end{defn}

\begin{defn}\label{def: 1-handle stabilization}
    Let $\Sigma$ and $W$ be as above, and let $\Sigma'$ be a $(g, n)$-stabilization of $\Sigma$. When $(g, n) = (0, 2)$ and $S_0 \cup D_1 \cup D_2$ bounds a 3-dimensional 1-handle embedded in $W$, we say $\Sigma'$ is a \emph{1-handle stabilization} of $\Sigma$.
\end{defn}

We will simply write ``stabilization'' instead of $(g, n)$-stabilization, and state explicitly when we mean 1-handle stabilization. We now formally define the two notions of stabilization distance. For simplicity, we will only define the stabilization distance for disks.

\begin{defn}\label{def: 1-handle stabilization distance}
    Let $\Sigma$ and $\Sigma'$ be disks in $W$ such that $\partial \Sigma = \partial \Sigma'$ and $[\Sigma] = [\Sigma'] \in H_2(W, \partial W; \Z)$. The \emph{1-handle stabilization distance} $d(\Sigma, \Sigma')$ between $\Sigma$ and $\Sigma'$ is the minimal number $k$ such that $\Sigma$ and $\Sigma'$ become isotopic rel boundary after each is stabilized $k$ times. 
\end{defn}

\begin{defn}\label{def: stabilization distance}
    Let $\Sigma$ and $\Sigma'$ be disks in $W$ such that $\partial \Sigma = \partial \Sigma'$ and $[\Sigma] = [\Sigma'] \in H_2(W, \partial W; \Z)$. The \emph{stabilization distance} $\mu(\Sigma, \Sigma')$ between $\Sigma$ and $\Sigma'$ is defined to be the minimum of 
    \[
    \max\{g(\Sigma_1), \hdots, g(\Sigma_k)\}
    \]
    over all sequences of connected surfaces from $\Sigma = \Sigma_1$ to $\Sigma' = \Sigma_k$ in $W$ such that $\partial \Sigma_i = K$ for all $i$ and $\Sigma_i$ and $\Sigma_{i+1}$ are related by a stabilization or a destabilization. 
\end{defn}

Note that necessarily $d(\Sigma, \Sigma') \ge \mu(\Sigma, \Sigma')$. However, as the next example illustrates, the two notions are distinct.

\begin{exam}\label{ex: d > mu}
    The knot $K = 9_{46}$ is shown in Figure \ref{fig: 9_46 isotopy 1}. $K$ bounds an obvious torus in $S^3$. Moreover, by compressing the two circles which generate the first homology of this torus, we obtain two slice disks $D$ and $D'$ with boundary $K$. Both disks can be described as banded unlinks with a single band \cite{swenton_calculus_2knots_4space,HKM_iso_surf} (Figure \ref{fig: 9_46 isotopy 1}). It is shown in \cite{MP_stab_dist} that $d(D, D') = 1$. This can be seen as follows. 

    \begin{figure}[h!]
        \includegraphics[scale = 0.4]{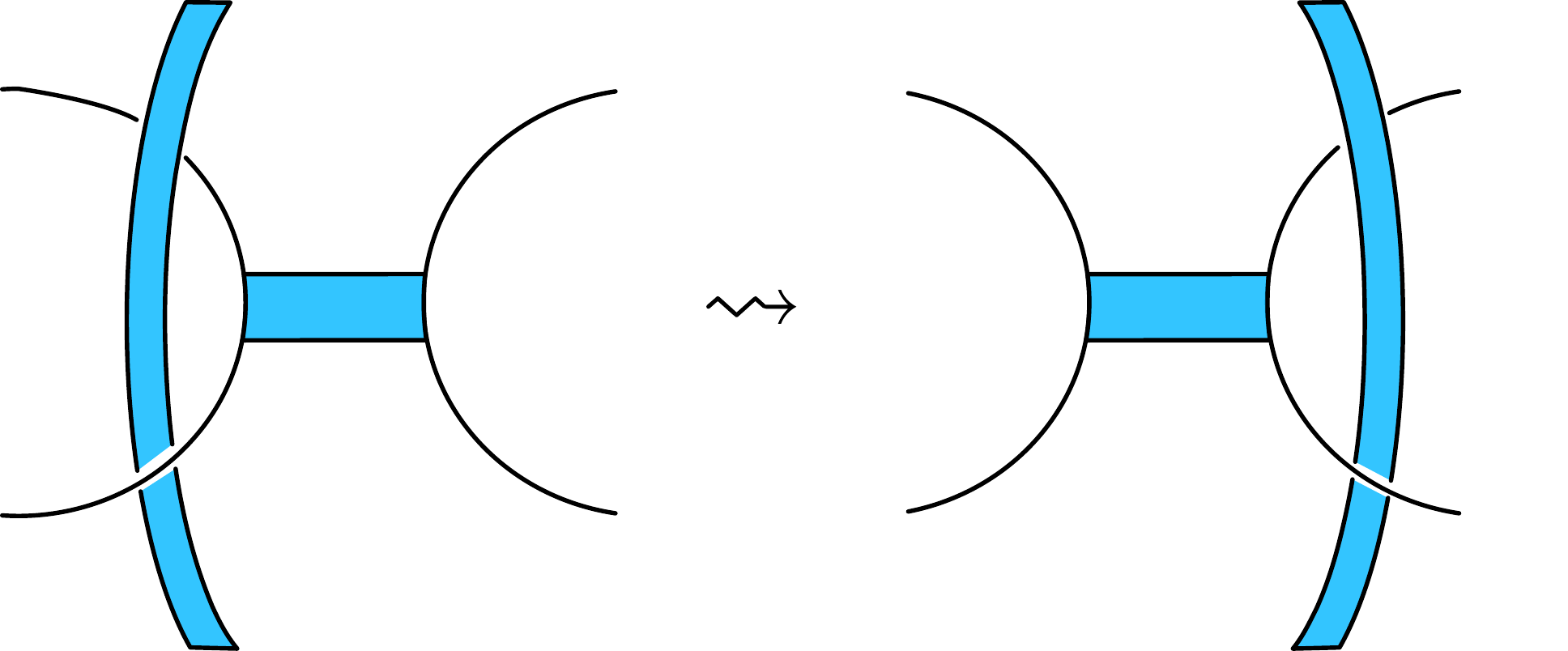}
        \caption{Swimming one band through another.}
        \label{fig: swim move}
    \end{figure} 

    \begin{figure}[h!]
        \includegraphics[scale=.7]{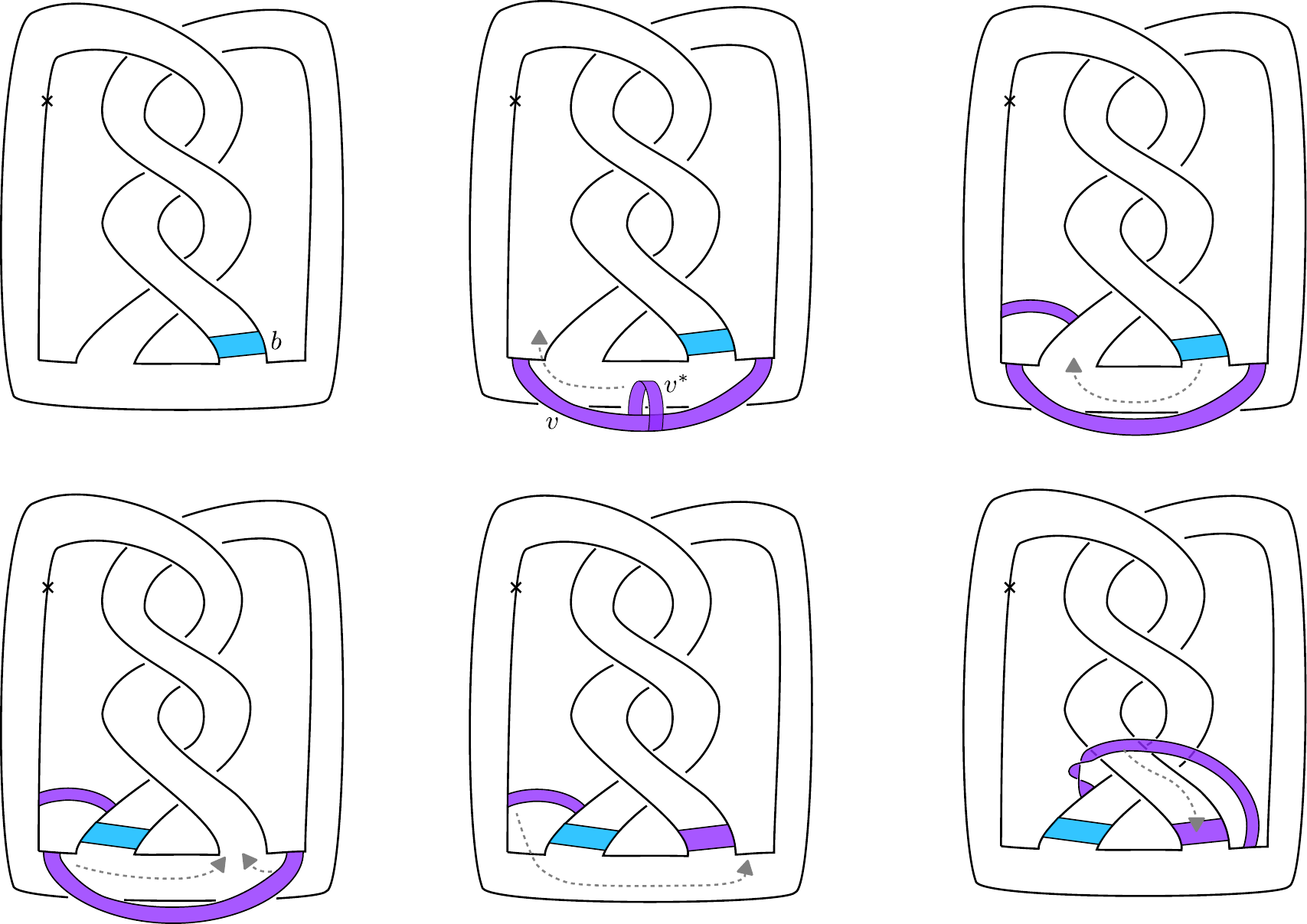}
        \caption{An explicit isotopy taking a 1-handle stabilization of $D$ to a 1-handle stabilization of $D'$. A swim move occurs in frame 6. (Continued in Figure \ref{fig: 9_46 isotopy 2}).}
        \label{fig: 9_46 isotopy 1}
    \end{figure}

    To show that the two disks become isotopic after a single 1-handle stabilization, it suffices to show that by attaching a tube, the relative positions of the two bands can be swapped, and then that the tube can be isotoped until it is once again clearly the result of a 1-handle stabilization.

    \begin{figure}[h!]
        \includegraphics[scale=.7]{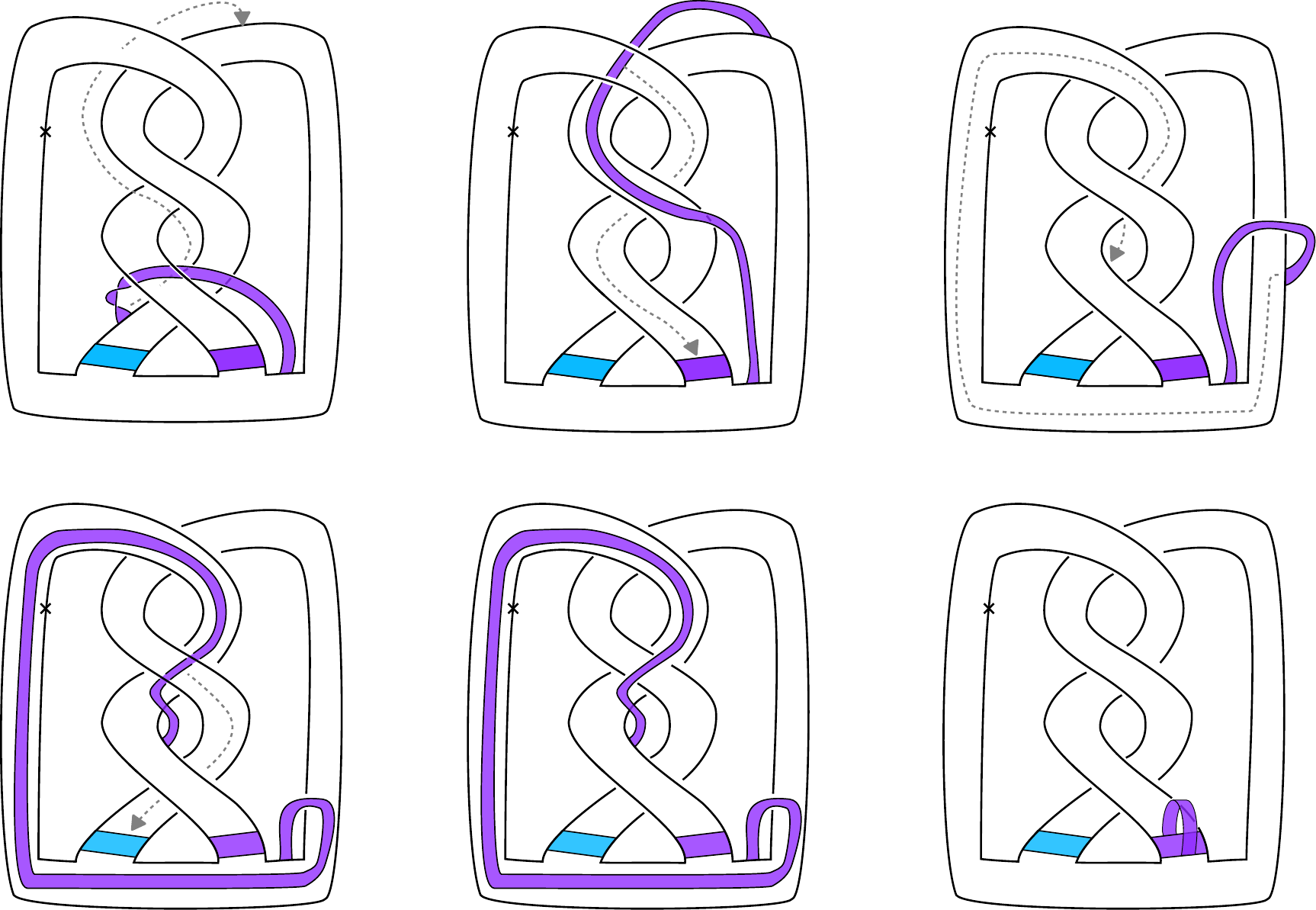}
        \caption{The remainder of the isotopy between the stabilizations of $D$ and $D'$. Swim moves occur in frames 2 and 4.}
        \label{fig: 9_46 isotopy 2}
    \end{figure}

    By strategically attaching a tube, we can slide the right-hand band, $b$, to the left (Figure \ref{fig: 9_46 isotopy 1}.) However, this band slide separates the tube into two bands, $v$ and $v^*$. Next, slide $v$ into the position originally occupied by $b$. Now, we can drag $v^*$ around $K$ until it once again forms a tube with $v$ by performing swim moves as necessary to change crossings of $v^*$ with the diagram for $K$. See Figure \ref{fig: swim move} for an example of a swim move, and see \ref{fig: 9_46 isotopy 2} for the remainder of the isotopy taking the stabilization of $D$ to the stabilization of $D'$. \\

    Miller-Powell use Alexander modules to show that by taking boundary connected sums of these disks, they can produce disks with arbitrarily large 1-handle stabilization distance: $d(\natural^m D,\natural^m D') = m$.  However, it is clear from Figures \ref{fig: 9_46 isotopy 1} and \ref{fig: 9_46 isotopy 2}, that $\mu(\natural^m D,\natural^m D') = 1$; no band ever slid over the marked point on the diagram, so we can take the connected sums at the marked points. Since the more general stabilization distance allows us to stabilize \emph{and} destabilize, we can attaching a single tube in order to isotope the first $D$ summand to $D'$, then destabilize, and then repeat the strategy on the next copy of $D$, until we are left with $\natural^m D'$. \\

    It is worth noting that even though the obvious ribbon disk for $\natural^m D$ has $m$ bands, we can in some sense ``reuse'' (stabilize then destabilize) the same tube $m$ times to make $\natural D$ and $\natural D'$ isotopic. $K$ has fusion number 1, so necessarily $\Ord_U(K) = 1$. By the knot Floer homology K\"unneth formula, 
    \[
        \Ord_U(\#^m K) =  \Ord_U(K) = 1. 
    \]
    This implies that $\tau(\natural^m D, \natural^m D') \le 1$ (see Definition \ref{def: secondary tau invt}), and therefore $\tau$ cannot detect the large $1$-handle stabilization distance of these disks. 
\end{exam}

\subsection{Concordances induced by cables}\label{subsection: prelims, satellite concordances}

Let $C: K \ra K'$ be a concordance. Given a pattern knot $P \sub S^1\times D^2$, we obtain a concordance between the satellites of $K$ and $K'$ as follows. Remove a neighborhood of $C$ in $S^3 \times I$ . The Seifert framing of $K$ determines an identification $\vphi:\partial \nu(K) \ra S^1\times \partial D^2$. Define the satellite concordance $C_P$ to be $(S^3\times I - C, \emptyset) \cup_{\vphi \times \id} (S^1\times D^2 \times I, P \times I).$ 

The $(p, q)$-cable of a knot is a satellite with pattern $P = T_{p, q}$, the $(p, q)$-torus knot. Since the $(p, 1)$-torus knot is the unknot, it is clear that if a knot $K$ is slice, so is its $(p, 1)$-cable. Moreover, by capping off the $(p, 1)$-cable of the unknot, we obtain a cabled disk for $K$. \\

By the work of Freedman and Quinn, locally flat proper submanifolds have topological normal bundles which are unique up to ambient isotopy \cite[Section 9.3]{freedman_quinn_top_4mflds}. Therefore, by the nature of the construction, topological isotopy is preserved by the cabling operation. Therefore, since the disks $D$ and $D'$ which bound the knot $J$ are topologically isotopic by the work of Conway and Powell \cite[Theorem 1.2]{conway_powell_characterisation_homotopy_ribbon}, this topological isotopy also produces a topological isotopy between the $(p, 1)$-cables of $D$ and $D'$. This gives us:

\begin{lemma} \label{lemma: cabled disks are topologically isotopic}
The cabled disks $D_p$ and $D_p'$ which bound $K_p$ are topologically isotopic. 
\end{lemma}

\subsection{An upper bound for the stabilization distance}\label{subsection: prelim, upper bound}

We now turn to the proof that the stabilization distance between $D_p$ and $D_p'$ is at most $p$.

As a warm up case, consider the disks $D$ and $D'$ with boundary $J$.  Let $b$ be the left-hand band which defines $D$ and let $b'$ be the right-hand band which defines $D'$. As in Example \ref{ex: d > mu}, attach a tube $v \cup v^*$ to $J$, so that the band $b$ can be slid until it becomes isotopic to the band $b'$. Next, slide the band $v$ into the position originally occupied by $b$. We have exchanged the roles of $b$ and $b'$ at the cost of tangling the bands $v$ and $v^*$ which made up the stabilization; it is not clear whether the resulting surface is isotopic to a stabilization of $D'$. If can isotope $v^*$ away from $J \cup b$ and back onto $v$, we will be done. 

\begin{figure}[h!]
    \centering
    \includegraphics[scale=.55]{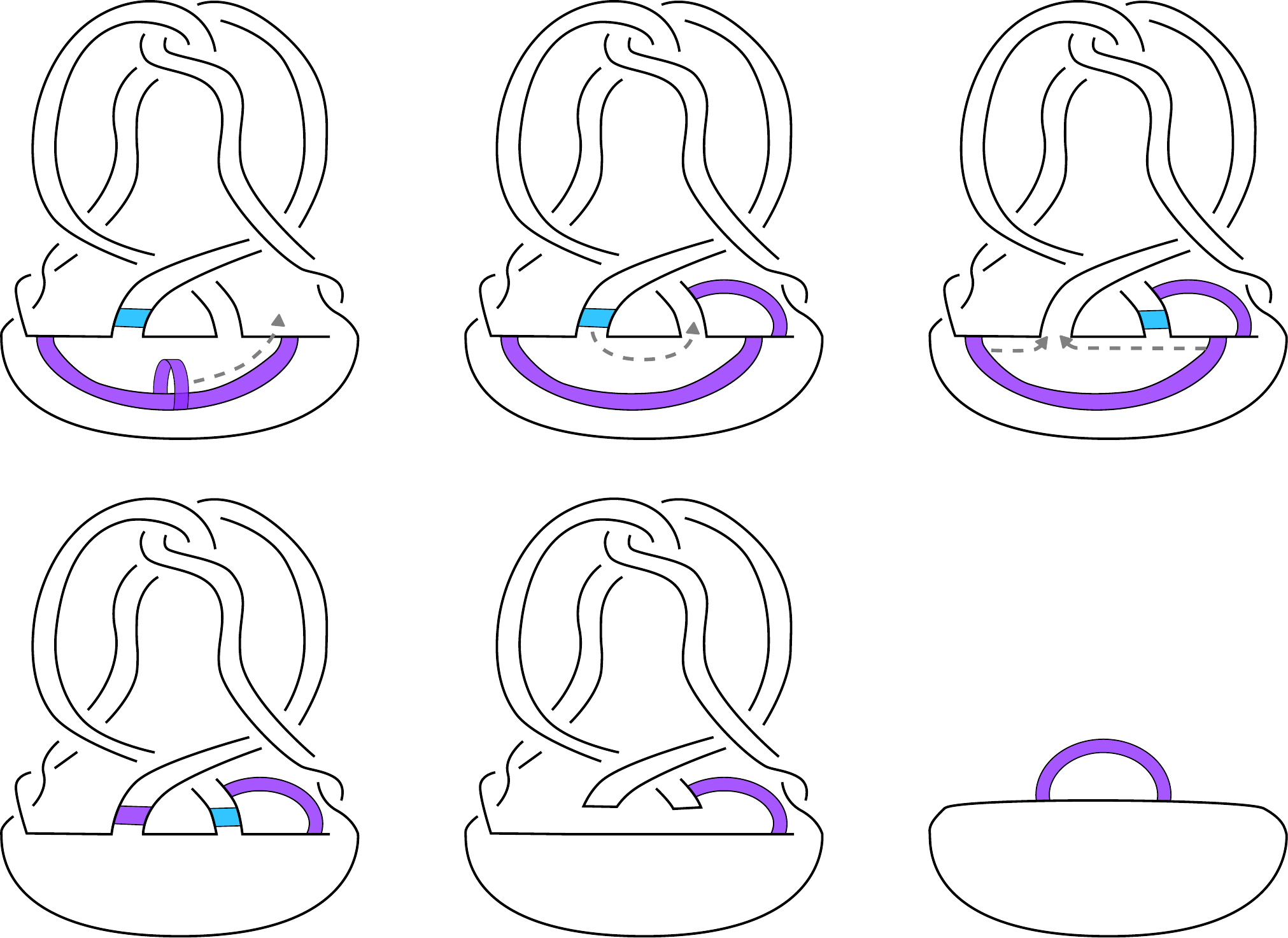}
    \caption{An isotopy of $D \cup (v \cup v^*)$.}
    \label{fig: warm up isotopy}
\end{figure}

As in Example \ref{ex: d > mu}, we can use the bands $b$ and $v$ to pull $v^*$ into the correct position by a sequence of swim moves. Recall that a swim move of $v^*$ through $b$ corresponds to pushing $v^*$ below the critical point for $b$, and performing an isotopy of $v^*$ in $S^3-J(b)$, i.e. in the complement of link obtained by band surgery on $b$. Performing the entire isotopy at this level will turn out to be easier to visualize, especially once we progress to the cabled case.

Push the band $v^*$ into the interior of $B^4$, (say to radius $r = 2/3$) below the critical points for the bands $b$ and $v$. Here, the level sets of the surface are isotopic to $J(b)(v) = J(b)(b')$ (the result of band surgery on both $b$ and $b'$) which is the unknot. Moreover, from Figure \ref{fig: warm up isotopy}, we see that the band $v^*$ is attached trivially.

At this point, the diagram is symmetric. Hence, this argument can be repeated with the disk $D'$ stabilized with tube $u \cup u^*$: we isotope $u^*$ into $B^4$ as well, until we see $u^*$ attached to the unknot. $U \cup v^*$ and $U\cup u^*$ are clearly isotopic, so by composing the first isotopy with the inverse of the second we obtain the desired isotopy from $D \cup (v\cup v^*)$ to $D' \cup (u \cup u^*)$.

The cabled case is similar.

\begin{figure}[h!]
    \includegraphics[scale=.65]{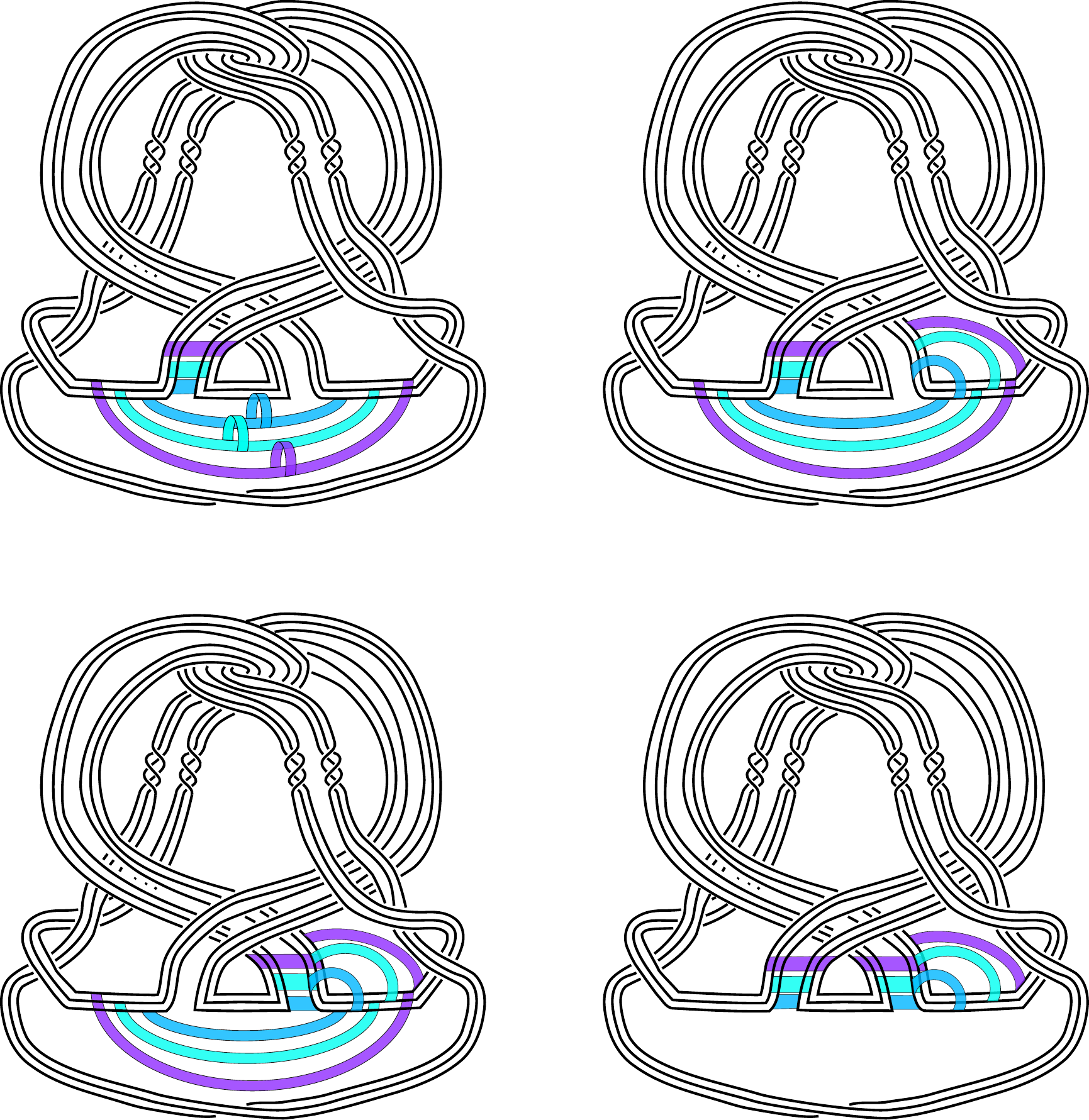}
    \caption{Part 1 of an isotopy between $p$-fold stabilizations of $D_p$ and $D_p'$}
    \label{fig: isotopy between stabilizated cabled disks 1} 
\end{figure}

\begin{figure}[h!]
    \includegraphics[scale=.65]{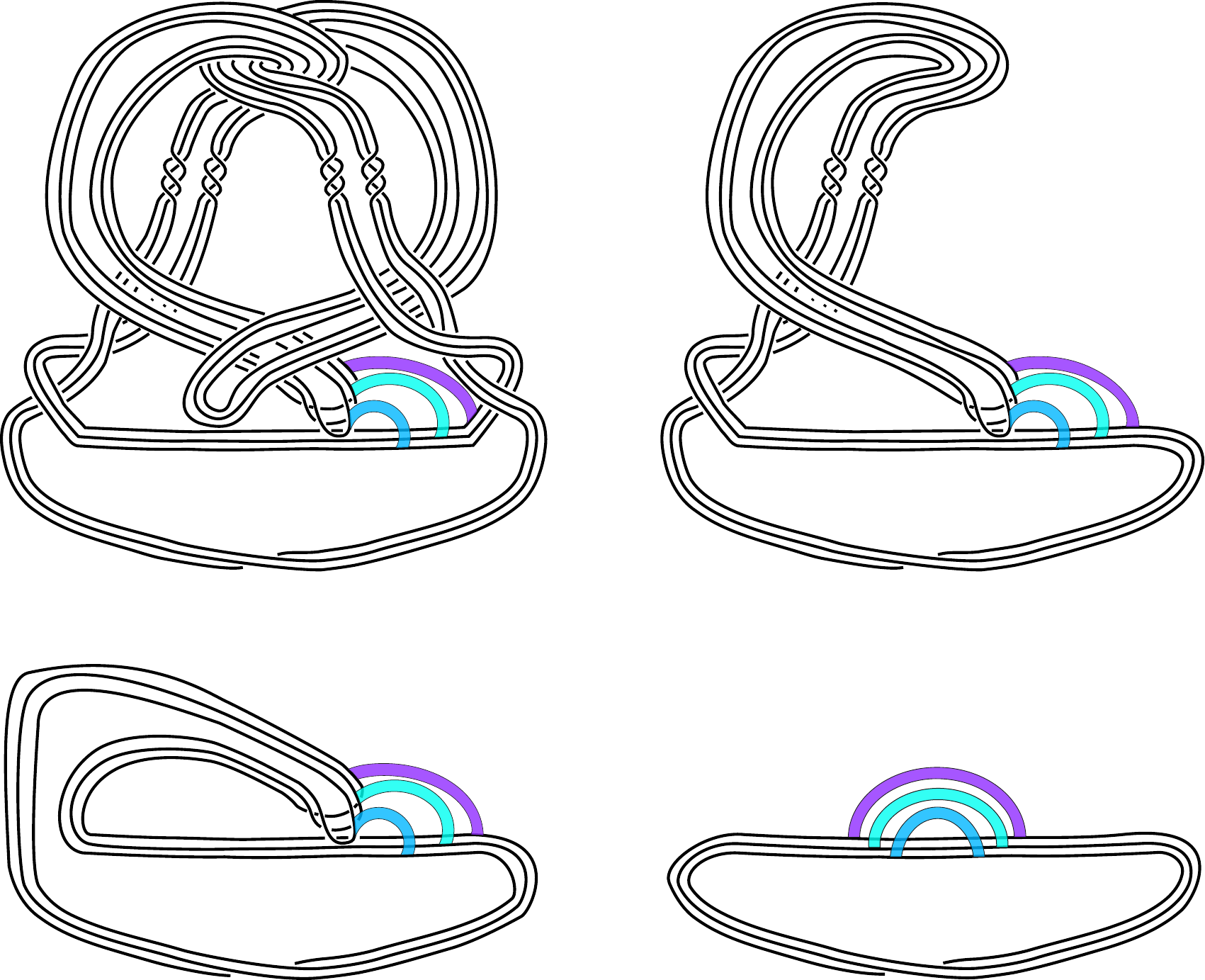}
    \caption{Part 2 of an isotopy between $p$-fold stabilizations of $D_p$ and $D_p'$}
    \label{fig: isotopy between stabilizated cabled disks 2} 
\end{figure}

\begin{prop}\label{prop: upper bound for cabled stabilization distance}
    The disks $D_p$ and $D_p'$ become isotopic after $p$ 1-handle stabilizations.
\end{prop}
\begin{proof}
    Figure \ref{fig: isotopy between stabilizated cabled disks 1} gives a band presentation for $D_p$. Let $b_1, \hdots, b_p$ and $b_1', \hdots, b_p'$ be the bands of $D_p$ and $D_p'$. As in the case $p = 1$, attach tubes $v_1 \cup v_1^*, \hdots, v_p \cup v_p^*$ in order to move the band $b_i$ into the position of $b_i'$, and then slide $v_i$ into the original position of $b_i$ (Figure \ref{fig: isotopy between stabilizated cabled disks 1}). The result of band surgery on all the $b_i$ and $v_i$ bands is again the unknot, and as shown in Figure \ref{fig: isotopy between stabilizated cabled disks 2} this unknot is naturally identified with the $(p, 1)$-cable of the unknot. Moreover, at this level set, the $v_i^*$ bands are attached trivially. So, just as before, perform the symmetric isotopy of $D_p'$, and concatenate to obtain an isotopy between the $p$-fold stabilizations of $D_p$ and $D_p'$. 
\end{proof}
\section{Floer-theoretic Background}

Theorem \ref{Theorem 1: exotic disks with large stabilization distance} is proved using a range of Floer-theoretic techniques. We begin this section by reviewing knot Floer homology, cobordism maps, and the Juh\'asz-Zemke stabilization distance bounds. Next, we review the techniques of equivariant knot Floer homology developed by Dai-Mallick-Stoffregen \cite{DaiMallickStoffregen_equivar_knots_floer} (which fit naturally within the context of $\iota$-complexes \cite{HMZ_connected_sum_involutive}). Finally, we review some bordered Floer homology and recall the Lipshitz-Ozsváth-Thurston pairing theorems \cite{LOT_bordered_HF}.

\subsection{Knot Floer homology}\label{subsection: background, knot floer}

Knot Floer homology is an invariant of knots in 3-manifolds defined by Ozsváth and Szabó \cite{os_knotinvts} and independently Rasmussen \cite{rasmussen_knotcompl}. The knot Floer groups are constructed by choosing a doubly-pointed Heegaard diagram $(\Sigma, \alpha, \beta, w, z)$ for $(Y, K)$, where $\alpha = (\al_1, ..., \al_{g})$ and $\beta = (\be_1, ..., \be_{g})$ are the attaching curves and $g$ is the genus of $\Sigma$. Denote by $\T_\al$ and $\T_\be$ the half dimensional tori $\al_1 \times ... \times \al_{g}$ and $\be_1 \times ... \times \be_{g}$ in $\Sym^{g}(\Sigma)$. Define $\cCFK^-(Y, \frs)$ to be the free $\F_2[U, V]$-module generated by intersection points $\x$ in $\T_\al\cap\T_\be$ with $\frs_w(\x) = \frs.$ The differential is defined by counting holomorphic disks of Maslov index 1, i.e.
\[
\partial(\x) = \sum_{\y\in \T_\al\cap \T_\be}\sum_{\substack{\phi \in \pi_2(\x, \y),\\ \mu(\phi) =1}} \# \widehat{\mathcal{M}}(\phi)U^{n_w(\phi)}V^{n_z(\phi)}\y,
\]

Maps between knot Floer complexes are induced by decorated link cobordisms.

\begin{defn}\label{Def: decorated link cob}
A \textit{decorated link cobordism} from $(Y_0, \bK_0) = (Y_0, (K_0, w_0, z_0))$ to $(Y_1, \bK_1) = (Y_1, (K_1, w_1, z_1))$ is a pair $(W, \cF) = (W, (\Sigma, \cA))$ with the following properties:
\begin{enumerate}
	\item $W$ is an oriented cobordism from $Y_0$ to $Y_1$
	\item $\Sigma$ is an oriented surface in $W$ with $\partial \Sigma = -K_0 \cup K_1$
	\item $\cA$ is a properly embedded 1-manifold in $\Sigma$, dividing it into subsurfaces $\Sigma_\w$ and $\Sigma_\z$ such that $w_0, w_1 \sub \Sigma_w$ and $z_0,z_1 \sub \Sigma_\z.$
\end{enumerate}
\end{defn}
In \cite{zemke_linkcob}, it is shown that a decorated link cobordism $(W, \cF)$ from $(Y_0, \bK_0)$ to $(Y_1, \bK_1)$ and a $\SpinC$-structure $\frs$ on $W$, give rise to a map
\[
F_{W, \cF, \frs}: \cCFK^-(Y_0, \bK_0, \frs|_{Y_0}) \ra \cCFK^-(Y_1, \bK_1, \frs|_{Y_1}),
\]
and these maps are functorial \cite[Theorem B]{zemke_linkcob}. The decorated link cobordism maps are defined as compositions of maps associated to handle attachments to the embedded surfaces and to the ambient 4-manifold. \\

We will make use of an algebraic variant of $\cCFK^-$, which is usually denoted $\CFK^-$. Define $\CFK^-(Y, K, \frs)$ be the $\F_2[U]$-module obtained from $\cCFK^-(Y, K, \frs)$ by setting $V = 0$ with differential 
\[
\partial(\x) = \sum_{\y\in \T_\al\cap \T_\be}\sum_{\substack{\phi \in \pi_2(\x, \y),\\ \mu(\phi) =1, \\ n_\z(\phi) = 0}} \# \widehat{\mathcal{M}}(\phi)U^{n_\w(\phi)}\y.
\]
Let $\HFK^-(Y, K, \frs)$ be the homology of this complex.

\subsection{Secondary invariants of a surface}

We now recall the definition of the secondary invariant $\tau$ of Juh\'asz and Zemke \cite{juhasz_zemke_stabilization_bounds}. Let $\Sigma$ be a surface in $B^4$ with boundary $\bK = (K, w, z)$, a doubly based knot. Decorate $\Sigma$ by a single arc such that the $\z$-subregion $\Sigma_\z \sub S$ is a bigon. Let $\cF = (\Sigma, \cA)$ be the resulting decorated surface. Then, there is an induced map 
\[
F_{B^4, \cF}: \F[U, V] \ra \cCFK^-(K, w, z).
\]

\begin{figure}[h!]\centering
	\includegraphics[scale=.5]{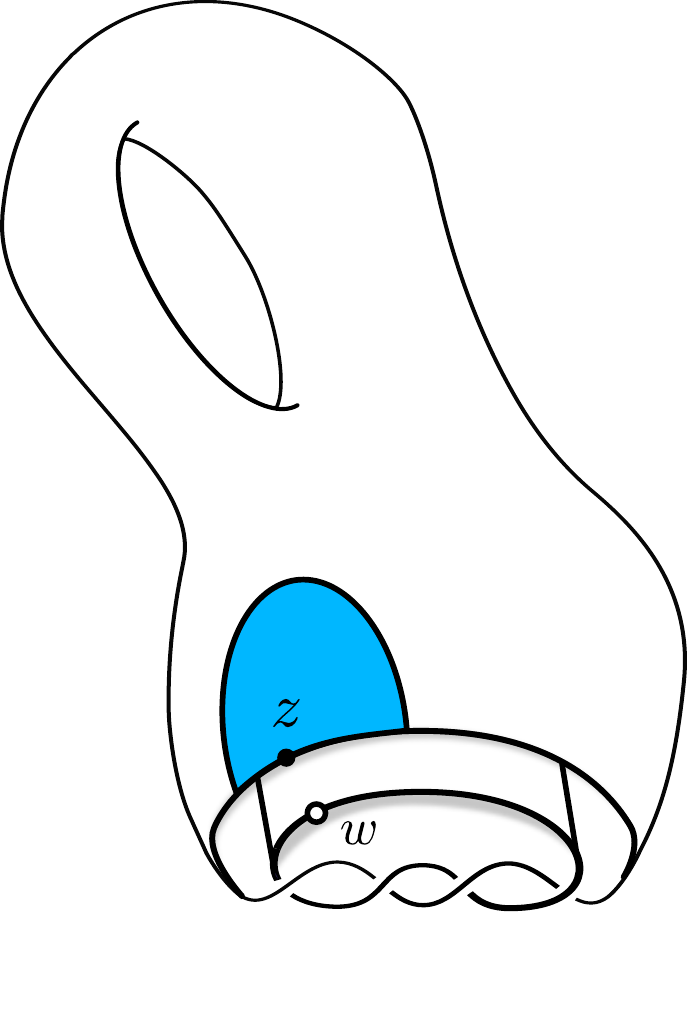}
	\caption{A surface $\Sigma$ with $\Sigma_\z$-sub region a bigon.}
	\label{fig: decorated surface}
\end{figure}

\begin{definition}\label{def: secondary tau invt}{\cite[Definition 4.4]{juhasz_zemke_stabilization_bounds}}
	Let $(K, w, z)$ be a doubly based knot in $S^3$ and let $\Sigma$ and $\Sigma'$ be two disks in $B^4$ with boundary $K$, decorated as above. Then, define
	\[
	\tau(\Sigma, \Sigma') = \min\{ n : U^{n} \cdot[F_{B^4, \cF}(1)] = U^{n}\cdot[F_{B^4, \cF'}(1)] \in \HFK^-(K) \}.
	\]
\end{definition}

\begin{remark}
	This invariant, $\tau$ (and its relatives $\nu, V_k, \Upsilon$), are called ``secondary invariants'' because they are defined as surface analogues of existing knot invariants.
\end{remark}

A key result of \cite{juhasz_zemke_stabilization_bounds} states this invariant provides a lower bound for the stabilization distance. 

\begin{theorem}{\cite[Theorem 1.1]{juhasz_zemke_stabilization_bounds}}
	Let $K$ be in a knot in $S^3$, and let $\Sigma, \Sigma'$ be disks in $B^4$ with boundary $K$. Then, 
	\[
	\tau(\Sigma, \Sigma') \le \mu(\Sigma, \Sigma'),
	\]
	where $\mu(\Sigma, \Sigma')$ is the stabilization distance between $S$ and $S'$.
\end{theorem}

\subsection{Equivariant Floer homology and $\iota$-complexes:}

\begin{defn}
An \emph{$\i$-complex}, is a chain complex $(C, \partial)$, which is free and finitely generated over $\F[U]$, together with an endomorphism $\i: C \ra C$ satisfying the following properties:
\begin{enumerate}
    \item $C$ has a $\Z$-grading, with respect to which $U$ has grading $-2$, which we call the Maslov grading.
    \item There is an isomorphism $U^{-1}H_*(C) \ra \F[U, U^{-1}]$.
    \item $\i$ preserves the grading, and, up to homotopy, squares to the identity: $\i^2 \simeq \id$.
\end{enumerate}
\end{defn}

Following Hendricks and Manolescu \cite{hendricks_manolescu_Invol}, numerical invariants can be extracted from an $\i$-complex in the following way. Define 
\[
I^-(C) = \Cone(C \xra{Q(1 + \i)} Q\cdot C),
\]
where $Q$ is a formal variable which squares to zero. Then define 
\[
\underline{d}^\i = \max\{r| \exists x \in H_*(I^-(C))_r, \forall n, U^n x \ne 0 \text{ and } U^n x \not \in \im(Q)\}-1
\]
\[
\overline{d}^\i = \max\{r| \exists x \in H_*(I^-(C))_r, \forall n, U^n x \ne 0 \text{ and } U^n x \in \im(Q)\}.
\]
These are called $\i$-correction terms, and are invariant under $\i$-homotopy equivalence.\\

In fact, the correction terms are invariant under the following, weaker form of equivalence.

\begin{defn}
Let $(C_1, \i_1)$ and $(C_2, \i_2)$ be $\i$-complexes. 
\begin{enumerate}
    \item A \emph{local map} from $(C_1, \i_1)$ to $(C_2, \i_2)$ is a grading preserving $\i$-homomorphism $f: C_1 \ra C_2$ which induces an isomorphism from $U^{-1}H_*(C_1)$ to $U^{-1}H_*(C_2)$.
    \item The $\i$-complexes $(C_1, \i_1)$ and $(C_2, \i_2)$ are called \emph{locally equivalent} if there is a local map from $(C_1, \i_1)$ to $(C_2, \i_2)$ as well as a local map from $(C_2, \i_2)$ to $(C_1, \i_1)$.
\end{enumerate}
\end{defn}

$\iota$-complexes arise naturally in Floer theory. Given a homology $3$-sphere $Y$, the $\SpinC$-conjugation action induces an action on $\CF^-(Y)$, which gives rise to an $\iota$-complex \cite{hendricks_manolescu_Invol}. For $N\ge g(K)$, the correction terms 
\[
\underline{V_0}^\iota(K) :=-\frac{1}{2}\underline{d}^\iota(S_N^3(K), [0]),\,\, \overline{V_0}^\iota(K) :=-\frac{1}{2}\overline{d}^\iota(S_N^3(K), [0]) 
\]
give bounds on the slice genus of $K$ \cite{juhász_zemke_newslicegenus}. \\

Similarly, given an involution $\tau$ on a homology 3-sphere $Y$, the action of $\tau$ on $CF^-(Y)$ again gives rise to an $\iota$-complex \cite{DaiHeddenMallick_corks_invol_floer}. As in the $\SpinC$-conjugation case, for $N\ge g(K)$, the correction terms 

\[
\underline{V_0}^\tau(K) :=-\frac{1}{2}\underline{d}^\tau(S_N^3(K), [0]),\,\, \overline{V_0}^\tau(K) :=-\frac{1}{2}\overline{d}^\tau(S_N^3(K), [0]) 
\]
give bounds on the equivariant slice genus of $K$. In fact, these numerical invariants bound a slightly more subtle quantity, which they call the \emph{isotopy-equivariant slice genus}, which is defined to be the minimal genus of a surface $\Sigma$ such that $\tau(\Sigma)$ is isotopic to $\Sigma$ rel boundary. They show knot Floer homology can distinguish the exotic disks $D$ and $D'$ by showing the isotopy-equivariant genus of $J$ is non-zero. \\

In analogy with the refinement of involutive Floer homology for knots, there is an algebraic notion of an $\ik$-complex. 

\begin{defn}
An $\ik$-complex $(C_K, \partial, \ik)$ is a free, finitely generated chain complex over $\F[U, V]$, equipped with an endomorphism $\ik$ satisfying the following properties:
\begin{enumerate}
    \item $\ik$ is skew-$\F[U, V]$-equivariant, i.e. $U\circ \ik = \ik \circ V$ and $V\circ \ik = \ik \circ U$.
    \item Up to homotopy, $\ik$ squares to the Sarkar map $\id + \Phi\Psi$.
\end{enumerate}
\end{defn}

The $\SpinC$-conjugation involution induces an action $\iota_K$ on $\cCFK^-(K)$ giving it the structure of an $\iota_K$-complex. Similarly, a geometric involution on $K$ induces an action $\tau_K$ on $\cCFK^-(K)$ as well, which is closely related to an $\iota_K$-complex (in the case $K$ is strongly invertible $\tau_K^2 \simeq \id$, not the Sarkar map).

Notice that $\ik$ commutes with $UV$, so an $\ik$-complex $(C, \partial, \ik)$ can be viewed as an infinitely generated complex with an $\F[U]$-equivariant endomorphism $\ik$ (where $U$ acts as $UV$). Equip $C$ with an Alexander and Maslov grading, and define $\bA_0^-(C)$ to be the subset of $C$ in Alexander grading zero. $\bA_0^-(C)$ is an $\F[U]$-module, and $(\bA_0^-(C), \iota_K)$ is an ordinary $\iota$-complex. Moreover, for the $\SpinC$-conjugation action, the action of a geometric involution, and their composition, there are large surgery formulas \cite{hendricks_manolescu_Invol,Mallick_surgery_equivariant_knots} \[(\bA_0^-(\cCFK^-(K)), \sigma_K) \simeq (\CF^-(S_N^3(K), [0]), \sigma), \text{ for } N \ge g(K),\] for $\sigma \in \{\iota, \tau,  \iota \circ \tau\}$. 

Going forward, we will denote $\bA_0^-(K) := \bA_0^-(\cCFK^-(K))$.

\subsection{Bordered Floer homology}

Bordered Floer homology is a package of invariants of 3-manifolds with parametrized boundary. Bordered Floer homology associates to a surface $F$ a differential graded algebra $\cA(F)$ and to a 3-manifold $Y$ with boundary together with an identification $\vphi: \partial Y \ra F$ a left differential graded module over $\cA(-F)$, called $\CFDh(Y)$ and a right $\cA_\infty$-module over $\cA(F)$ called $\CFAh(Y)$. Much like classical Heegaard Floer homology, the bordered Floer invariants are defined by representing a 3-manifold with parametrized boundary by a kind of Heegaard diagram and counting holomorphic disks which, in the bordered case, may asymptotically approach the boundary.\\

Bordered Floer homology has a pairing theorem \cite[Theorem 1.3]{LOT_bordered_HF}, which recovers the hat-version of the Heegaard Floer homology of the manifold obtained by gluing bordered manifolds along their common boundary. Given 3-manifolds $Y_1$ and $Y_2$ with $\partial Y_1 \cong F \cong \partial Y_2$, there is a homotopy equivalence
\[
	\CFh(Y_1\cup Y_2) \simeq \CFAh(Y_1) \boxtimes_{\cA(F)} \CFDh(Y_2).
\]

Bordered Floer theory also recovers knot Floer homology \cite[Theorem 11.21]{LOT_bordered_HF}. Given a doubly pointed bordered Heegaard diagram $(\cH_1, w, z)$ for $(Y_1, \partial F, K)$ and a bordered Heegaard diagram $(\cH_2, z)$ with $\partial Y_1 \cong F \cong -\partial Y_2$, then 
\[
	\HFK^-(Y_1 \cup Y_2, K) \cong H_*(\CFA^-(\cH_1, w, z) \boxtimes_{\cA(F)} \CFDh(\cH_2, z)).
\]
Bordered Floer theory, therefore, gives an effective way to study satellites. Let $K_P$ be the satellite of $K$ with pattern $P$. $\CFDh(S^3-K)$ is determined by $\CFK^-(K)$ \cite[Chapter 11]{LOT_bordered_HF}, $\HFK^-(K_P)$ can be computed by finding a doubly pointed bordered Heegaard diagram $\cH_P$ for the pattern $P$ in the solid torus and computing the box tensor product $\CFA^-(\cH_P) \boxtimes \CFDh(S^3-K)$.\\

The last important result from bordered Floer theory which we will utilize, is the morphism spaces pairing theorem, which gives a means of recovering classical Heegaard Floer homology in terms of the $\Hom$ functor rather than the tensor product functor \cite{LOT_HF_as_morphism}. If $Y_1$ and $Y_2$ are 3-manifolds with $\partial Y_1 \cong F \cong \partial Y_2$, then 
\[
\CFh(-Y_1 \cup Y_2) \simeq \Mor_{\cA(-F)}(\CFDh(Y_1), \CFDh(Y_2)), 
\]
where the latter object is the chain complex of $\cA(-F)$-linear maps from $\CFDh(Y_1)$ to $\CFDh(Y_2)$ equipped with differential 
\[
d(\vphi) = \partial_{\CFDh(Y_2)} \circ \vphi + \vphi \circ \partial_{\CFDh(Y_1)}.
\]

\section{Maps associated to exotic disks}\label{Section: computations of exotic disk maps}

In this section, we use the strong involution on the knot $J$ to study the maps $F_D$ and $F_{D'}$. 

\begin{prop}\label{prop: compute F_D and F_D'}
	Let $D$ and $D'$ be the exotic pair of disks bound by $J$. Then, 
	\[
		F_D(1) + F_{D'}(1)= e_1 + e_2.
	\]
\end{prop}

\subsection{Computing the disk maps on $\CFK^-(K)$:}

The exotic disks $D$ and $D'$ are relatively complicated, so computing these maps is nontrivial. The key topological property of these disks, which makes computing these maps feasible, is their symmetry, i.e. the existence of a strong involution $\tau$ on $K$ which extends to a involution on $B^4$ which takes the disk $D$ to $D'$. 

To compute the disk maps, we use the fact that 
\begin{equation}\label{eqn: naturality of involution}
	[F_{D'}(1)] = [F_{\tau D}(1)] =  \tau_K([F_{D}(1)]) \in H_*(\cCFK^-(J)).
\end{equation}
See \cite[Lemma 7.9]{DaiMallickStoffregen_equivar_knots_floer}. Dai-Mallick-Stoffregen show that 
\[\min\{ \underline{d}^\tau(S_1^3(J), [0]), \underline{d}^{\iota\tau}(S_1^3(J), [0]) \} \le -2, \] which, by definition, implies 
\begin{equation}\label{eqn: nontrivial V invt}
	\max\{\underline{V}_0^{\tau}(J), \underline{V}_0^{\iota\tau}(J)\} \ge 1.
\end{equation}
Since this quantity provides a lower bound for the isotopy equivariant slice genus, this implies there is no smooth isotopy of $D$ taking it to $D'$.

\begin{figure}
	\centering
	\scalebox{.8}{\begin{tikzpicture}
        \node at (2, 2) (e2) {$e_2$};
        \node at (1.5, 1.5) (e1) {$e_1$};
        \node at (1, 1) (x) {x};
        \node at (.5, .5) (a1) {$a_1$};
        \node at (0, 0) (a2) {$a_2$};
        \node at (0-3, .5) (Ub1) {$Ub_1$};
        \node at (.5-3, 0) (Ub2) {$Ub_2$};
        \node at (.5, 0-3) (Vc1) {$\;\;Vc_1$};
        \node at (0, .5-3) (Vc2) {$Vc_2\;\;$};
        \node at (-3, 0-3) (UVe1) {$\;\;UVe_1$};
        \node at (.5-3, .5-3) (UVe2) {$\;\;UVe_2$};
        \node at (.5-4, .5-4) (UVx) {$UV x$};
        \node[scale=1.25] at (-4, -4) (dots) {$\iddots$};

        \draw[thick, ->] (a1) -- (Ub1);
        \draw[thick, ->] (a1) -- (Vc1);
        \draw[thick, ->] (a2) -- (Ub2);
        \draw[thick, ->] (a2) -- (Vc2);

        \draw[thick, ->] (Ub1) -- (UVe1);
        \draw[thick, ->] (Ub2) -- (UVe2);
        \draw[thick, ->] (Vc1) -- (UVe1);
        \draw[thick, ->] (Vc2) -- (UVe2);

    \renewcommand{\d}{8}
        \node at (1.5+\d, 1.5) (j2) {$j_2$};
        \node at (1+\d, 1) (j1) {$j_1$};
        \node at (.5+\d, .5) (f1) {$f_1$};
        \node at (0+\d, 0) (f2) {$f_2$};
        \node at (0-3+\d, .5) (Ug1) {$Ug_1$};
        \node at (.5-3+\d, 0) (Ug2) {$Ug_2$};
        \node at (.5+\d, 0-3) (Vh1) {$\;\;Vh_1$};
        \node at (0+\d, .5-3) (Vh2) {$Vh_2\;\;$};
        \node at (-3+\d, 0-3) (UVj1) {$\;\;UVj_1$};
        \node at (.5-3+\d, .5-3) (UVj2) {$\;\;UVj_2$};
        \node[scale=1.25] at (-3.5+\d, -3.5) (dots) {$\iddots$};

        \draw[thick, ->] (f1) -- (Ug1);
        \draw[thick, ->] (f1) -- (Vh1);
        \draw[thick, ->] (f2) -- (Ug2);
        \draw[thick, ->] (f2) -- (Vh2);

        \draw[thick, ->] (Ug1) -- (UVj1);
        \draw[thick, ->] (Ug2) -- (UVj2);
        \draw[thick, ->] (Vh1) -- (UVj1);
        \draw[thick, ->] (Vh2) -- (UVj2);

        \node at(2.75, -1.2) {$\bigoplus$};

        \foreach \j in {0-1, 1-1, 2-1, 3-1}{

        \node at (-3 + 4*\j, -6) {$\bullet$};
        \node at (-3 + 4*\j, -5.5) {$e_1$};
        \node at (-2 + 4*\j, -6) {$\bullet$};
        \node at (-2 + 4*\j, -5.5) {$x$};
        \node at (-1 + 4*\j, -6) {$\bullet$};
        \node at (-1 + 4*\j, -5.5) {$e_2$};
        \node at (-3 + 4*\j, -7) {$\bullet$};
        \node at (-3 + 4*\j, -6.5) {$j_1$};
        \node at (-1 + 4*\j, -7) {$\bullet$};
        \node at (-1 + 4*\j, -6.5) {$j_2$};

        \draw[very thick, -] (-2 + 4*\j, -6) -- (-2 + 4*\j, -9);
        }


        \node at (-2 + 4*0-1, -4.5) {$1$};
        \node at (-2 + 4*1-1, -4.5) {$\iota$};
        \node at (-2 + 4*2-1, -4.5) {$\sigma$};
        \node at (-2 + 4*3-1, -4.5) {$\iota \circ \sigma$};


        \draw[thick, blue!80, ->] (-2.2 + 4*0-1, -5.25) to[out=100, in=80,looseness=3] (-1.8 + 4*0-1, -5.25);
        \draw[thick, blue!80, ->] (-3.2 + 4*0-1, -5.25) to[out=100, in=80,looseness=3] (-2.8 + 4*0-1, -5.25);
        \draw[thick, blue!80, ->] (-1.2 + 4*0-1, -5.25) to[out=100, in=80,looseness=3] (-.8 + 4*0-1, -5.25);


        \draw[thick, blue!80, <->] (-3 + 4*1-1, -5.25) to [bend left=90] (-1 + 4*1-1, -5.25);
        \draw[thick, blue!80, ->] (-2.2 + 4*1-1, -5.25) to[out=100, in=80,looseness=3] (-1.8 + 4*1-1, -5.25);


        \draw[thick, blue!80, <->] (-3 + 4*2-1, -5.25) to [bend left=90] (-1 + 4*2-1, -5.25);
        \draw[thick, blue!80, ->] (-2.2 + 4*2-1, -5.25) to[out=100, in=80,looseness=3] (-1.8 + 4*2-1, -5.25);
        \draw[thick, blue!80, ->] (-2.2 + 4*2-1, -5.5) -- (-2.8 + 4*2-1, -5.5);
        \draw[thick, blue!80, ->] (-1.8 + 4*2-1, -5.5) -- (-1.2 + 4*2-1, -5.5);


        \draw[thick, blue!80, ->] (-1.2 + 4*3-1, -5.25) to[out=100, in=80,looseness=3] (-.8 + 4*3-1, -5.25);
        \draw[thick, blue!80, ->] (-2.2 + 4*3-1, -5.25) to[out=100, in=80,looseness=3] (-1.8 + 4*3-1, -5.25);
        \draw[thick, blue!80, ->] (-3.2 + 4*3-1, -5.25) to[out=100, in=80,looseness=3] (-2.8 + 4*3-1, -5.25);
        \draw[thick, blue!80, ->] (-2.2 + 4*3-1, -5.5) -- (-2.8 + 4*3-1, -5.5);
        \draw[thick, blue!80, ->] (-1.8 + 4*3-1, -5.5) -- (-1.2 + 4*3-1, -5.5);
\end{tikzpicture}}
	\caption{The complex $\bA_0^-(J)$ as well all possible involutions on its homology. The solid lines indicate the action of $U$.}\label{fig: tau action}
\end{figure}
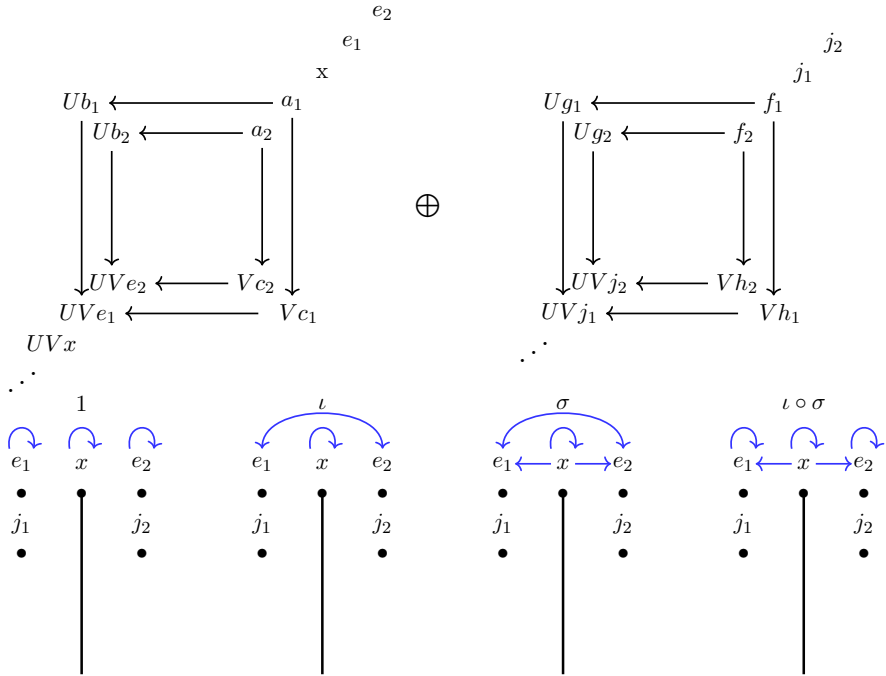

But, they deduce more. The knot Floer complex can of $J$ be computed using Szab\'o's knot Floer calculator \cite{szabo_knot_floer_calc}, and is shown in Figure \ref{figure:CFK(J)}. For more details, see the proof of \cite[Theorem 1.6]{DaiMallickStoffregen_equivar_knots_floer}. With this information, $\bA_0^-(J)$ can be computed, as well as its homology. See Figure \ref{fig: tau action}.

\begin{figure}[h!]
\begin{minipage}{.4\textwidth}
\[
x \quad \ \begin{tikzcd}[row sep=1cm, column sep=1cm] b_i\ar[d,"V"] & a_i \ar[l, "U"] \ar[d, "V"]\\
e_i & c_i \ar[l, "U"] \\
 g_i\ar[d,"V"] & f_i \ar[l, "U"] \ar[d, "V"]\\
j_i & h_i \ar[l, "U"] 
\end{tikzcd}
\]
\end{minipage}
\begin{minipage}{.2\textwidth}
\[
\begin{tabular}{ c | c | c }
\text{Generator} & $\gr_U$ & $\gr_V$ \\
\hline
$x$ & 0 & 0 \\
\hline
$a_i$ & 0 & 0 \\
$b_i$ & $1$ & $-1$ \\
$c_i$ & $-1$ & $1$ \\
$e_i$ & 0 & 0 \\
\hline
$f_i$ & $-1$ & $-1$ \\
$g_i$ & 0 & $-2$ \\
$h_i$ & $-2$ & 0 \\
$j_i$ & $-1$ & $-1$
\end{tabular}
\]
\end{minipage}

\caption{On the left is the complex $\cCFK^-(J)$, $i \in \{1, 2\}$. On the right are the bi-gradings of generators.}\label{figure:CFK(J)}
\end{figure}
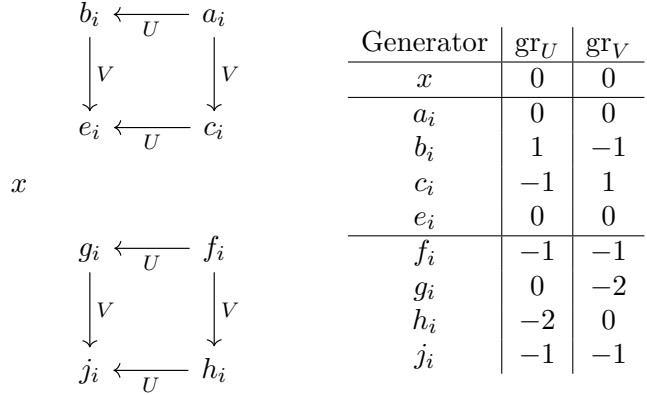

\begin{proof}[Proof of Proposition \ref{prop: compute F_D and F_D'}:]

The action of $\tau_K$ on the summand $A\sub \bA_0^-(J)$ generated by $x, a_i, b_i, c_i$, and $e_i$ can be determined by simply considering all possible involutions, and computing the associated correction terms. It is not hard to see that the map quotienting out the complement of this complex is a local equivalence. Therefore, since $\underline{V}_0^\tau$ and $\underline{V}_0^{\iota\tau}$ are local invariants, these numerical invariants will be unaffected by simply assuming $\tau_K$ acts as the identity on $j_i$.  

There are four possible involutions on $H_*(\bA_0^-(J))$ which must be considered, which are shown in Figure \ref{fig: tau action} as well. Clearly the first map is the identity. The second map is the $\SpinC$-conjugation action; even though this complex is not thin, it decomposes into a sum of thin complexes, so $\iota_K$ can be computed using the techniques of \cite[Section 8.3]{hendricks_manolescu_Invol}. Since $\underline{V}_0^{\id}$ and $\underline{V}_0^{\iota}$ are concordance invariants, the sliceness of $J$ implies $\underline{V}_0^{\id}(J)$ and $\underline{V}_0^{\iota}(J)$ are both zero. A direct computation shows $\underline{V}_0^{\iota \sigma} (J) = 0$ as well and that $\underline{V}_0^{\sigma} (J) = 1$. Therefore, by Equation \ref{eqn: nontrivial V invt}, the involution $\sigma$ corresponds to \emph{either} $\tau_K$ or $\iota\circ \tau_K$, or equivalently, $\tau_K$ is equal to $\sigma$ or $\iota \circ \sigma$. \\

Suppose $\tau_K = \sigma$. By \cite[Corollary 1.8]{zemke_absgr}, the map induced by doubling either disk is the identity. Therefore, $F_D(1)$ and $F_{D'}(1)$ are of the form $x + \vep_1 \cdot e_1 + \vep_2 \cdot e_2$ where $\vep_i \in \{0, 1\}$. By \cite[Theorem 1.6]{DaiMallickStoffregen_equivar_knots_floer}, $F_D(1)$ and $F_{D'}(1)$ are distinct, and by Equation \ref{eqn: naturality of involution}, are exchanged by the action of $\tau_K$. Since the involution $\sigma$ acts as the identity on the elements $x + e_1$ and $x + e_2$, the only possibility is that $\{F_D(1), F_{D'}(1)\} = \{x, x + e_1 + e_2\}$.

Now, suppose $\tau_K = \iota \circ \sigma$. By the same argument, $F_D(1)$ and $F_{D'}(1)$ are of the form $x + \vep_1 \cdot e_1 + \vep_2 \cdot e_2$ where $\vep_i \in \{0, 1\}$. But, in this case, there are two possibilities, since under the action of $\iota\circ \sigma$, $x + e_1$ and $x + e_2$ are exchanged, as are $x$ and $x + e_1 + e_2$. Therefore, either $\{F_D(1), F_{D'}(1)\} = \{x, x + e_1 + e_2\}$ or $\{F_D(1), F_{D'}(1)\} = \{x + e_1, x + e_2\}$.

Therefore, in both cases, the sum of the two maps is $e_1 + e_2$.
\end{proof}

\section{Satellite Concordance Maps}
In this section, we prove Theorem \ref{Theorem 2: satellite formula for concordance maps}, which will allow us to compute the maps associated to satellite concordances. Since it requires no additional effort, we prove Theorem \ref{Theorem 2: satellite formula for concordance maps} for homology concordances, i.e. link cobordisms $(W, C): (Y, K) \ra (Y', K')$ where $Y$ and $Y'$ are integer homology spheres, $C$ is an annulus, and $W$ is a homology cobordism. 

Let $(W, C): (Y, K) \ra (Y', K')$ be a homology concordance. Let $\cF = (C, \cA)$ be the annulus $C$ decorated with a pair of parallel arcs running from $K$ to $K'$. This decorated concordance induces a map $\CFK^-(Y, K) \ra \CFK^-(Y', K')$, denoted $F_{W,C}$, which is defined as a composition of elementary cobordism maps \cite{zemke_linkcob}. Choose a Morse function $f: W \ra \R$ on $W$ so that $f|_C$ has no critical points and choose a gradient-like vector field which is tangent to $C$. This induces a handle decomposition for $(W, C)$ which only involves four-dimensional handles. 

Since the restriction of this Morse function to the surface has no critical points, this also produces a handle decomposition for the complement of $C$; the attaching curves for the handles are already embedded in the complement of $K$, so the cobordism $W - C$ is built simply by attaching the same handles to $Y-K$. 

\begin{prop}\label{Prop: concordance complement map}
    Let $(W, C): (Y, K) \ra (Y', K')$ be a homology concordance. Then, given a handle decomposition for $(W, C)$ as above, there exists a map $F_{W - C}$ with the property that for any pattern knot $P$ in the solid torus, the following diagram commutes up to homotopy:
    \begin{center}
    \begin{tikzcd}
        \CFA^-(\cH_P) \boxtimes \CFDh(Y-K) \ar[d, "\bI_{\cH_P} \boxtimes F_{W - C}"]  \ar[r, "\simeq"]  
            & \CFK^-(Y, K_P) \ar[d, "F_{W, C_P}"]  \\
        \CFA^-(\cH_P) \boxtimes \CFDh(Y'-K') \ar[r, "\simeq"]
            & \CFK^-(Y', K'_P).
    \end{tikzcd}    
    \end{center}
    Here, $\cH_P$ is a doubly pointed Heegaard diagram for the knot $P$ in the solid torus, $K_P$ and $K_P'$ are the satellites of $K$ and $K'$ with pattern $P$, $(W, C_P)$ is the concordance induced by the pattern, and the horizontal homotopy equivalences are given by the pairing theorem for knot Floer homology \cite{LOT_bordered_HF}.
\end{prop} 

The map $F_{W-C}: \CFDh(Y-K) \ra \CFDh(Y'-K')$ will be defined in the standard way, as a composition of maps corresponding to handle attachments \cite{os_holotri}.

\subsection{1- and 3-handle maps:}

The maps associated to 1- and 3-handle attachments are simplest to define. For simplicity, let $Y$ be an integer homology sphere. We will write $Y(S^0)$ for the result of $S^0$-surgery on $Y$ (which is, of course, diffeomorphic to $Y \# S^1\times S^2$).

\begin{lemma}\label{lemma: 1-handle map}
    Let $F: \CFK^-(Y, K_P)\ra\CFK^-(Y(S^0), K_P', \frt_0)$ be the 1-handle cobordism map. There exists a map $\widetilde{F}: \CFDh(Y-K) \ra \CFDh((Y-K)(S^0), \frt_0|_{(Y-K)(S^0)})$ making the following diagram commute up to homotopy:
    \begin{center}
    \begin{tikzcd}
        \CFA^-(\cH_P) \boxtimes \CFDh(Y-K) \ar[d, "\bI_{\cH_P} \boxtimes \widetilde{F}"]  \ar[r, "\simeq"]  
            & \CFK^-(Y, K_P) \ar[d, "F"]  \\
        \CFA^-(\cH_P) \boxtimes \CFDh((Y-K)(S^0), \frt_0|_{(Y-K)(S^0)}) \ar[r, "\simeq"]
            & \CFK^-(Y(S^0), K_P, \frt_0),
    \end{tikzcd}
    \end{center}
    where $\frt_0$ is the torsion $\SpinC$-structure on $Y(S^0)$, i.e. $c_1(\frt_0) =0$.
\end{lemma}

\begin{proof}
    Fix a nice Heegaard diagram $\cH = (\Sigma_g, \alpha_1^c, \hdots, \alpha_{g-1}^c, \beta_1, \hdots, \beta_g)$ for $Y-K$. Choose curves $\lambda$ and $\mu$ in $\Sigma$ isotopic to a longitude and meridian of $K$ respectively, which intersect in a single point and avoid the $\alpha$ circles. A bordered Heegaard diagram $\cH_{B}$ for $Y-K$ is obtained from $\cH$ by deleting a neighborhood of $p$, and defining $\alpha_1^a =  \lambda - p$ and $\alpha_2^a = \mu-p$ to be the $\alpha$-arcs parametrizing the boundary, i.e.
    \[
        \cH_{B} = (\Sigma - p, \alpha_1^a, \alpha_2^a, \alpha_1^c, \hdots, \alpha_{g-1}^c, \beta_1, \hdots, \beta_g).
    \]
    Let $\cH_P$ be a nice doubly pointed bordered Heegaard diagram for the pattern knot embedded in the solid torus. A doubly pointed Heegaard diagram $\cH_{K_P}$ for $(Y, K_P)$ is obtained by gluing $\cH_B$ and $\cH_P$ along their common boundary. \\

    Recall the definition of the 1-handle map $F: \CFK^-(Y, K_P) \ra \CFK^-(Y(S^0), K_P', \frt_0)$. Choose a pair of points $p_1$ and $p_2$ in $\cH_{K_p}$ away from $\cH_P$. Moreover, assume $p_1$ and $p_2$ lie in the same connected component of $\Sigma - (\cup_i \alpha_i) - (\cup_i \beta_i)$ as the basepoint $z$. Remove neighborhoods of $p_1$ and $p_2$, and attach an annulus. Add two new curves $\alpha_0$ and $\beta_0$ which are homologically essential in the annulus and intersect transversely in a pair of points, which we denote $\theta^+$ and $\theta^-$. There are two bigons from $\theta^+$ to $\theta^-$. The 1-handle map is simply 
    \[
        \x \mapsto \x \otimes \theta^+.
    \]
    In exactly the same way, $S^0$-surgery on the bordered Heegaard diagram for $Y-K$ gives rise to a map $\widetilde{F}: \CFDh(Y- K) \ra \CFDh((Y-K)(S^0)), \frt_0|_{(Y-K)(S^0)})$ defined $\x' \mapsto \x' \otimes \theta^+$. Since we chose nice diagrams, the identification of $\CFA^-(\cH_P) \boxtimes \CFDh(Y-K) \xra{\simeq} \CFK^-(Y, K)$ is simply the map which takes a pair of intersection points in $\cH_P$ and $\cH_K$ and views them as a single intersection point in $\cH_P \cup \cH_K$. Tautologically, then, the desired diagram commutes. 
\end{proof}
    
    The case of the 3-handles is dual to this case, and follows similarly. 

\begin{lemma}\label{lemma: 3-handle map}
    Let $H: \CFK^-(Y(S^0), K_P, \frt_0)\ra\CFK^-(Y, K_P')$ be the 3-handle cobordism map. There exists a map $\widetilde{H}: \CFDh((Y-K)(S^0), \frt_0|_{(Y-K)(S^0)}) \ra \CFDh(Y-K')$ making the following diagram commute up to homotopy:
    \begin{center}
    \begin{tikzcd}
        \CFA^-(\cH_P) \boxtimes \CFDh((Y-K)(S^0), \frt_0|_{(Y-K)(S^0)}) \ar[d, "\bI_{\cH_P} \boxtimes \widetilde{H}"]  \ar[r, "\simeq"]  
            & \CFK^-(Y(S^0), K_P, \frt_0) \ar[d, "H"]  \\
        \CFA^-(\cH_P) \boxtimes \CFDh(Y-K') \ar[r, "\simeq"]
            & \CFK^-(Y, K_P').
    \end{tikzcd}
    \end{center}
\end{lemma}

\subsection{2-handle map:}

The 2-handle cobordism is the only interesting case, and follows from the pairing theorem for triangles \cite[Proposition 5.35]{LOT_spectral_seq_II}. Before we define the 2-handle cobordism maps on $\CFDh$, we recall some facts about bordered Heegaard triple diagrams.

Given a doubly pointed bordered Heegaard diagram $\cH_{\alpha, \beta^0}$, we obtain a doubly pointed bordered Heegaard triple diagram $\cH_{\alpha, \beta^0, \beta^1}$ by performing a Hamiltonian translation on each of the $\beta^0$-curves. By removing the $\beta^0$-curves, we obtain an ordinary doubly pointed bordered Heegaard diagram $\cH_{\alpha, \beta^1}$. Counting holomorphic triangles defines a map
\[
m_2: \CFA^-(\cH_{\alpha, \beta^0})\otimes \CFA^-(\cH_{\beta^0, \beta^1}) \ra \CFA^-(\cH_{\alpha, \beta^1}).
\]
Taking $\Theta_{\beta^0, \beta^1}$ to be the top graded generator of the homology of $\CFh(\cH_{\beta^0, \beta^1})$, we can consider the map 
\[
m_2(-, \Theta_{\beta^0, \beta^1}): \CFA^-(\cH_{\alpha, \beta^0}) \ra \CFA^-(\cH_{\alpha, \beta^1}).
\]
We will make use of the fact that this map is homotopic to the map $\Psi_{\cH_{\alpha, \beta^1} \leftarrow \cH_{\alpha, \beta^0}}$ induced by the isotopy of $\beta^0$ to $\beta^1$, and is just the ``nearest point map'', taking an intersection point in $\alpha \cap \beta^0$ to the closest intersection point in $\alpha \cap \beta^1$. For a proof in the classical case, see \cite[Proposition 11.4]{Lipshitz_cylindrical} or \cite[Lemma 9.7]{JTZ_naturality_mapping_class_groups}. 

\begin{figure}[h!]\centering
    \includegraphics[scale=.5]{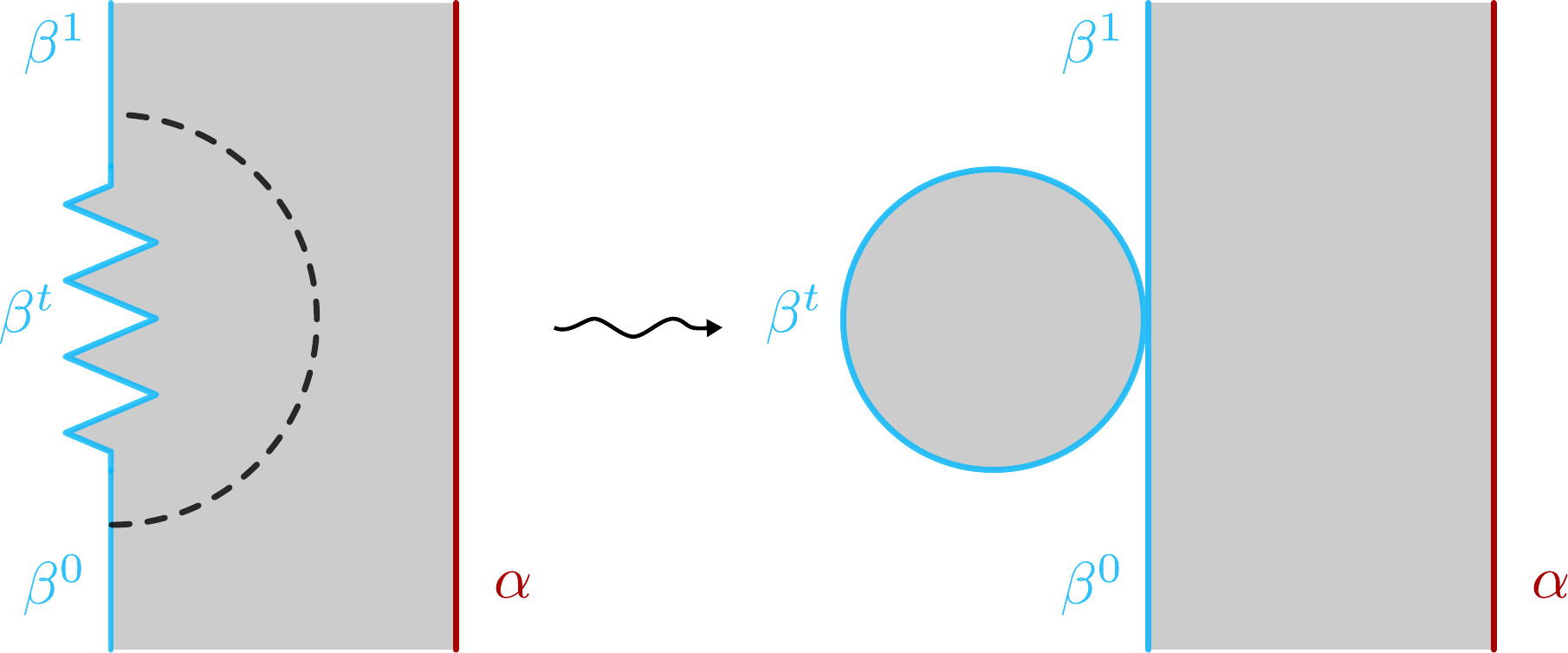}
    \caption{By pinching along the dotted line, we see a dynamic bigon map is homotopic to the composition of a monogon map with a triangle map.}
    \label{fig: neck stretching}
\end{figure}

\begin{lemma}\label{lemma: closest point map}
    Let $\cH_{\alpha, \beta^0, \beta^1}$ be a bordered Heegaard triple diagram where the $\beta^1$-curves are small Hamiltonian translates of the $\beta^0$-curves. Let $\Theta_{\beta^0, \beta^1}$ represent the top graded generator of $H_*(\CFh(\beta^0, \beta^1))$. Then, 
    \[
        m_2(-, \Theta_{\beta^0, \beta^1}) \sim \Psi_{\cH_{\alpha, \beta^1} \leftarrow \cH_{\alpha, \beta^0}}.
    \]
    Moreover, $\Psi_{\cH_{\alpha, \beta^1} \leftarrow \cH_{\alpha, \beta^0}}$ is the nearest point map.
\end{lemma}
\begin{proof}
    The holomorphic polygon maps are defined by counting certain holomorphic maps 
    \[
    u: (S, \partial S) \ra (\Sigma \times \Delta, (\alpha\times e_0) \cup (\beta^1 \times e_1)\cup\hdots \cup (\beta^n\times e_n), 
    \]
    where $\Delta$ is a disk with $n$ boundary punctures, and edges labeled $e_0, \hdots, e_n$. If $B \in \pi_2(x_0, \hdots, x_n; \rho_1, \hdots, \rho_m)$, then denote by $\cM^B(x_0, \hdots, x_n; \rho_{1}, \hdots, \rho_m)$ the moduli space of holomorphic maps in the homotopy class $B$ with asymptotics $x_1, \hdots, x_n, \rho_1, \hdots, \rho_m$. The map $\Psi_{\cH_{\alpha, \beta^1} \leftarrow \cH_{\alpha, \beta^0}}$ is defined by counting bigons with dynamic boundary conditions, i.e. maps
    \[
    (S, \partial S) \ra (\Sigma \times [0, 1] \times \R, (\alpha \times 1 \times \R) \cup \bigcup_t(\beta^t \times 0 \times \{t\}))
    \]
    where $\beta^{t} = \beta^0$ for $t \le 0$ and $\beta^{t} = \beta^1$ for $t \ge 1$. By a neck stretching argument, the moduli space of such maps splits into a product of $\cM^B(x, y, z; \rho_{1}, \hdots, \rho_m)$ and $\cM^A(x; \rho_{1}, \hdots, \rho_m)$. In other words, $\Psi_{\cH_{\alpha, \beta^1} \leftarrow \cH_{\alpha, \beta^0}}$ is homotopic to the composition of the monogon map,
    \[
    \theta: \F \ra \CFA^-(\cH_{\beta^0, \beta^1})
    \]
    with the triangle map $m_2(-,-)$. See Figure \ref{fig: neck stretching}. The monogon map simply picks out an element of $\CFA^-(\cH_{\beta^0, \beta^1})$. Therefore, $\Psi_{\cH_{\alpha, \beta^1}} \sim m_2(-, \theta(1)).$

    \textbf{Claim}: $\theta(1) = \Theta_{\beta^0, \beta^1}$. In particular, $m_2(-, \theta(1)) = m_2(-, \Theta_{\beta^0, \beta^1})$, as desired.\\

    If $\theta(1)$ was not equal to $\Theta_{\beta^0, \beta^1}$, then take the standard genus $g$ Heegaard diagram for $S^1\times S^2$ and consider the two maps \[m_2(-, \theta(1)), \Psi_{\cH_{\alpha, \beta^1} \leftarrow \cH_{\alpha, \beta^0}}: \CFA^-(\cH_{\alpha, \beta^0}) \ra \CFA^-(\cH_{\alpha, \beta^0})\] from above. The map $\Psi_{\cH_{\alpha, \beta^1} \leftarrow \cH_{\alpha, \beta^0}}$ is a homotopy equivalence, but is homotopic to $m_2(-, \theta(1))$. If $\theta(1)$ were not the highest graded generator, the composition would be zero, a contradiction.\\

    Finally, as $t \ra 0$, and $\beta^t$ approaches $\beta^0$, holomorphic disks with Maslov index 0 limit to bigons for the pair $(\alpha, \beta^0)$. But, any bigon for $(\alpha, \beta^0)$ with Maslov index 0 must be constant, since the $\R$-action on moduli space of homotopy classes of non-constant bigons is free, implying the dimension of the moduli space is nonzero.
\end{proof}

We are now ready to define the 2-handle cobordism maps. Let $\cL$ be a framed link in $Y$. Let $W(\cL)$ be the cobordism corresponding to attaching 2-handles along $\cL$. Let $Y(\cL)$ be the 3-manifold obtained by surgery on $\cL$. As usual, the 2-handle map is defined by counting holomorphic triangles. 

\begin{lemma}\label{lemma: 2-handle map}
    Fix a $\SpinC$-structure $\frs$ on $W(\cL)$. Let $G_\frs: \CFK^-(Y, K_P, \frs|_{Y})\ra\CFK^-(Y', K'_P, \frs|_{Y'})$ be the 2-handle cobordism map. There exists a map $\widetilde{G}_\frs: \CFDh(Y-K, \frs|_{Y-K}) \ra \CFDh(Y'-K', \frs|_{Y'-K'})$ making the following diagram commute up to homotopy:
    \begin{center}
    \begin{tikzcd}
        \CFA^-(\cH_p) \boxtimes \CFDh(Y-K, \frs|_{Y-K}) \ar[d, "\bI_{\CFA^-(\cH_P)} \boxtimes \widetilde{G}_\frs"]  \ar[r, "\simeq"]  
            & \CFK^-(Y, K_P, \frs|_{Y}) \ar[d, "G_\frs"]  \\
        \CFA^-(\cH_P) \boxtimes \CFDh(Y'-K', \frs|_{Y'-K'}) \ar[r, "\simeq"]
            & \CFK^-(Y', K'_P, \frs|_{Y'}).
    \end{tikzcd}
    \end{center}
\end{lemma}

\begin{proof}
    To define the map for 2-handles, let $\cL$ be a framed, $k$-component link in $Y - K$, and let $\cB$ be a bouquet for $\cL$. Choose a Heegaard triple diagram $\cH^{\alpha, \beta, \beta'}$ which is subordinate to this bouquet in the sense of \cite{os_holotri} and also admits a decomposition into two bordered Heegaard diagrams $\cH_P^{\alpha, \beta, \beta'} \cup \cH_B^{\alpha, \beta, \beta'}$, where $\cH_P^{\alpha, \beta, \beta'}$ is the bordered Heegaard triple diagram obtained from $\cH_P$ by adding Hamiltonian translates of the $\beta$-curves. Such a diagram can be produced as follows.
    \begin{enumerate}
        \item Start with a Heegaard diagram $\cH = (\Sigma_g, \alpha_1, \hdots, \alpha_{g-1}, \beta_{k+1}, \hdots, \beta_{g})$ for $Y - K -\cB$. For each $\beta$-curve in $\cH$, add a $\beta'$-curve which is a Hamiltonian translate so $|\beta_i \cap \beta'_j| = 2\delta_{ij}$.
        \item Add a collection $\{\beta_1, ..., \beta_k\}$ which are meridians of the components of $\cL$. Attaching 1- and 2-handles along the $\alpha$ and $\beta$ curves yields $Y-K$.
        \item For each component of $\cL$, add a curve $\beta_i'$ according to its framing. Attaching 1- and 2-handles along the $\alpha$ and $\beta'$ curves produces $(Y-K)(\cL) = Y'-K'$.
        \item Choose a longitude and meridian of $K$ in $\Sigma$ disjoint from the $\alpha$-curves which intersect in a single point, $p$. Define arcs $\alpha_1^a =  \lambda - p$ and $\alpha_2^a = \mu-p$.
    \end{enumerate}
    Altogether, this defines a bordered Heegaard triple diagram \[\cH_B^{\alpha, \beta, \beta'} = (\Sigma_g - p,\alpha_1^a, \alpha_2^a, \alpha_1^c, \hdots, \alpha_{g-1}^c, \beta_1, \hdots, \beta_{g}, \beta'_1, \hdots, \beta'_{g}).\] By construction, $\cH_P^{\alpha, \beta, \beta'} \cup \cH_B^{\alpha, \beta, \beta'}$ is a Heegaard triple subordinate to our chosen bouquet $\cB$.\\

    For $\delta, \vep \in \{\alpha, \beta, \beta'\}$, let $\cH_P^{\delta, \vep}$ and $\cH_B^{\delta, \vep}$ be the various standard bordered Heegaard diagrams associated to these triples. Let $\Theta_{\beta, \beta'}$ and $\Theta$ be the top dimensional generators of $H_*(\CFh(\cH_B^{\beta, \beta'}))$ and $H_*(\CFh(\cH_P^{\beta, \beta'}))$ respectively. The pairing theorem for triangles \cite[Proposition 5.35]{LOT_spectral_seq_II} gives the following homotopy-commutative square:
    \begin{center}
        \begin{tikzcd}
            \CFA^-(\cH_P^{\alpha, \beta}) \boxtimes \CFDh(\cH_B^{\alpha, \beta}) \ar[d, "m_2 (-\text{,} \Theta) \boxtimes m_2 (-\text{,} \Theta_{\beta, \beta'})" left]  \ar[r, "\simeq"]  
                & \CFK^-(\cH_P^{\alpha, \beta}\cup\cH_B^{\alpha, \beta}) \ar[d, "m_2 (-\text{,} \Theta\otimes\Theta_{\beta, \beta'})"]  \\
            \CFA^-(\cH_P^{\alpha, \beta'}) \boxtimes \CFDh(\cH_B^{\alpha, \beta'}) \ar[r, "\simeq"]
                & \CFK^-(\cH_P^{\alpha, \beta'}\cup\cH_B^{\alpha, \beta'}).
        \end{tikzcd}
    \end{center}

    The right vertical arrow is by definition the 2-handle cobordism map $G$ on $\CFK^-$. After identifying $\CFA^-(\cH_\infty^{\alpha, \beta})$ with $\CFA^-(\cH_\infty^{\alpha, \beta'})$, the map $m_2 (-, \Theta)$ is the identity. Define $m_2 (-, \Theta_{\beta, \beta'})$ to be the 2-handle cobordism map $\widetilde{G_\frs}$ on $\CFDh$. 
\end{proof}
    
Since each handle attaching cobordism map was defined with respect to a particular Heegaard diagram, the last step is to ensure maps induced by the Heegaard moves relating two diagrams induce homotopy equivalences which are compatible with the bordered Floer homology pairing theorem. This is shown by Hendricks-Lipshitz in \cite{HL_Inv_bordered_floer}.

\begin{lemma}{\cite[Lemma 5.6]{HL_Inv_bordered_floer}}\label{lemma: change diagram map}
    Suppose $\cH_1$ and $\cH_2$ are a pair of bordered Heegaard diagrams related by a bordered Heegaard move and $\cH_0$ is another bordered Heegaard diagram with $\partial \cH_0 = - \partial \cH_i$, $i \in \{1, 2\}$. Then, the diagram 
\end{lemma} 
    \begin{center}
        \begin{tikzcd}
            \CFA^-(\cH_1) \boxtimes \CFDh(\cH_0) \ar[d]  \ar[r, "\simeq"]  
                & \CFK^-(\cH_1\cup \cH_0) \ar[d]  \\
            \CFA^-(\cH_2) \boxtimes \CFDh(\cH_0) \ar[r, "\simeq"]
                & \CFK^-(\cH_2\cup \cH_0)
        \end{tikzcd}
    \end{center}
    commutes up to homotopy. The vertical maps come from the proof of invariance of bordered and classical Floer homology.

\begin{proof}[Proof of Proposition \ref{Prop: concordance complement map}:]
    We have a decomposition of $W - C$ into handle-attachment cobordisms, $W_1 \cup W_2 \cup W_3$. Define $F_{WC}$ to be the composition 
    \[
        F_{W-C} = \widetilde{H}\circ \Psi_{\cH_3 \leftarrow \cH_2} \circ \widetilde{G}_\frs\circ \Psi_{\cH_2 \leftarrow \cH_1}\circ \widetilde{F},
    \]
    where $\Psi_{\cH_{i+1} \leftarrow \cH_i}$ $i \in \{1, 2\}$ are the change of diagram homotopy equivalences and $\frs$ is the restriction of the unique $\SpinC$-structure on $W$ to $W_2$. By stacking the diagrams from Lemma \ref{lemma: 1-handle map}, \ref{lemma: 3-handle map}, \ref{lemma: 2-handle map}, and \ref{lemma: change diagram map} the result follows. 
\end{proof}
\section{Computations from Bordered Floer Homology}\label{Section: Bordered Floer Computations}

We now turn to the task of computing the maps $F_{D_p}$ and $F_{D_p'}$ which are induced by the $(p, 1)$-cables of the exotic disks $D$ and $D'$. It is more convenient to view these disks as concordances $C_p$ and $C_p'$ from the unknot. $F_{C_p}$ and $F_{C_p'}$ will be computed in terms of $F_{S^3\times I - C}$ and $F_{S^3\times I - C'}$, which as we will see, are determined by $F_{C}$ and $F_{C'}$ (up to some indeterminacy). 

As stated in Section \ref{Section: computations of exotic disk maps}, $\cCFK^-(J)$ consists of a singleton generator $x$ as well as four boxes. The summands generated by the two boxes generated by $\{f_i, g_i, h_i, j_i\}$ contain no elements of bigrading $(0, 0)$, and therefore do not intersect the images of the maps $F_C$ and $F_{C'}$. For this reason, we will work primarily with the subcomplex of $\cCFK^-(J)$ generated by $\{x, a_i, b_i, c_i, e_i\}$ to simplify the notation.\\

\cite{LOT_bordered_HF} gives an algorithm for determining $\CFDh(S^3-K)$ in terms of $\cCFK^-(K)$. Figure \ref{figure: CFK to CFD} shows how a box and a singleton generator give rise to summands of $\CFDh(S^3-J)$.

\begin{figure}[h!]
\begin{minipage}{.2\textwidth}
\[
\begin{tikzcd} b_i\ar[d,"V"] & a_i \ar[l, "U"] \ar[d, "V"]\\
e_i & c_i \ar[l, "U"] 
\end{tikzcd} \quad x \;\; \rightsquigarrow
\]
\end{minipage}
\begin{minipage}{.1\textwidth}
\hspace{2cm}
\end{minipage}
\begin{minipage}{.5\textwidth}
\[
\begin{tikzcd}[row sep=1.1cm, column sep=1.1cm] b_i \ar[d, "\rho_1"]
            & y^1_i \ar[l, "\rho_2"]
                & a_i \ar[l, "\rho_3"] \ar[d, "\rho_1"]& \\
        y^2_i 
            &     
                & y^4_i 
                    & x  \ar[loop right, "\rho_{12}"] \\
        e_i \ar[u, "\rho_{123}"]
            & y^3_i \ar[l, "\rho_2"]
                & c_i \ar[l, "\rho_3"] \ar[u, "\rho_{123}"]
                    &
\end{tikzcd}
\]
\end{minipage}
\caption{On the left is the summand of $\cCFK^-(J)$ containing $F_{D}(1)$ and $F_{D'}(1)$. On the right is a model for the corresponding summand of $\CFDh(S^3 - J)$.}\label{figure: CFK to CFD}
\end{figure}
%

\subsection{Computing the morphism complex}\label{subsection: bordered, morphism}

By Theorem \ref{Theorem 2: satellite formula for concordance maps}, there exists a map $F: \CFDh(S^3-U) \ra \CFDh(S^3-J)$ with the property that for any pattern $P$ in the solid torus, $\bI_{\CFA^-(P)} \boxtimes F$ computes the concordance map induced by the pattern. In particular, if we take $P$ to be the longitudinal unknot in the solid torus, which we denote $(T_\infty, \lambda)$ this map also computes $F_C$ (respectively $F_{C'}$). Since we did not show this map $F$ is unique, we will try to pin it down by computing the morphism space $\Mor_{\cA(T^2)}(\CFDh(S^3-U), \CFDh(S^3-J))$ and considering which maps $f \in \Mor_{\cA(T^2)}(\CFDh(S^3-U), \CFDh(S^3-J))$ have the property that $\bI_{\CFA^-(T_\infty, \lambda)}\boxtimes f \simeq F_C$ (respectively $F_{C'}$).

We begin by computing the dimension of the homology of the morphism space $\Mor_{\cA(T^2)}(\CFDh(S^3-U), \CFDh(S^3-J)).$ 

\begin{lemma}\label{lemma: dimension count}
    The space of homotopy classes of maps from $\CFDh(S^3-U)$ to $\CFDh(S^3-J)$ is 10 dimensional.
\end{lemma}
\begin{proof}
    By \cite{LOT_HF_as_morphism}, there is a homotopy equivalence:
    \begin{align*}
        \Mor_{\cA(T^2)}(\CFDh(S^3-U), \CFDh(S^3-J)) &\simeq \CFh(-(S^3-U) \cup (S^3-J)) \\
        & = \CFh(S_0^3(J)).
    \end{align*}
    The mapping cone formula \cite{os_HFK_integer_surgeries} shows that $\HF^+(S^3_0(J), [0])$ is the homology of the mapping cone $H_*(\Cone(\bA_0^+ \xra{v_0 + h_0} \bB^+_0))$. We illustrate part of the complex. 

\begin{center}
\scalebox{0.7}{%
\begin{tikzpicture}
        \node at (1.5, 1.5 + 3) (b1) {$b_i$};
        \node at (1.5+3, 1.5) (c1) {$c_i$};
        \node at (1.5+3, 3.5+1) (UVa1) {$(UV)^{-1}a_i$};
        \node at (2+3, 3+2) (UVx) {$(UV)^{-1}x_i$};
        \node at (1.5, 1.5) (e1) {$e_i$};
        \node at (1, 1) (x) {$x_i$};
        \node at (.5, .5) (a1) {$a_i$};
        \node at (.5-3, .5) (Ub1) {$Ub_i$};
        \node at (.5, .5-3) (Vc1) {$\;\;Vc_i$};
        \draw[thick, ->] (a1) -- (Ub1);
        \draw[thick, ->] (a1) -- (Vc1);
        \draw[thick, ->] (UVa1) -- (b1);
        \draw[thick, ->] (UVa1) -- (c1);
        \draw[thick, ->] (b1) -- (e1);
        \draw[thick, ->] (c1) -- (e1);

    \renewcommand{\d}{11}
        \node[scale=1.5] at (-3, 7) (A0) {$\bA_0^+$};
        \node[scale=1.5] at (9, 7) (B0) {$\bB_0^+$};

        \node at (1.5 + \d, 1.5 + 3) (b1t) {$\tilde{b}_i$};
        \node at (1.5+3. + \d, 1.5) (c1t) {$\tilde{c}_i$};
        \node at (1+3.5+ \d, 3.5+1) (UVa1t) {$(UV)^{-1}\tilde{a}_i$};
        \node at (1+4+ \d, 4+1) (UVxt) {$(UV)^{-1}\tilde{x_i}$};
        \node at (1.5+ \d, 1.5) (e1t) {$\tilde{e}_i$};
        \node at (1+ \d, 1) (xt) {$\tilde{x_i}$};
        \node at (.5+ \d, .5) (a1t) {$\tilde{a}_i$};
        \node at (.5+ \d, .5-3) (Vc1t) {$\;\;V\tilde{c}_i$};
        \draw[thick, ->] (a1t) -- (Vc1t);
        \draw[thick, ->] (UVa1t) -- (b1t);
        \draw[thick, ->] (UVa1t) -- (c1t);
        \draw[thick, ->] (b1t) -- (e1t);
        \draw[thick, ->] (c1t) -- (e1t);

        \draw[thick, ->] (6, 0) -- (9, 0) node[midway, above] {$v_0 + h_0$};
\end{tikzpicture}
}
\end{center}

    The homology of the portion of the complex shown is $\cT^+ \langle x_i \rangle \oplus \cT^+ \langle \tilde{x}_i \rangle \oplus \F\langle Ub_i = Vc_i \rangle$, where  $\cT^+ = \F[U, U^{-1}]/(U \cdot \F[U])$. The homology of the summand which is not shown (the remaining two boxes) is $\F\langle Ug_i = V f_i \rangle$. From $\HF^+(S_0^3(J))$, $\HFh(S^3_0(J))$ is obtained via the exact triangle 
    \begin{center}
        \begin{tikzcd}[column sep = small]
        \HF^+(S_0^3(J)) \ar[rr, "U"]  & & 
            \HF^+(S_0^3(J)) \ar[dl] \\
            & \HFh(S_0^3(J)). \ar[ul]
    \end{tikzcd}
    \end{center}
    
    A straightforward computation shows that $\HFh(S^3_0(J)) = \ker(U) \oplus \text{coker}(U) \cong \F^{\oplus 10}$. 
\end{proof}

Having computed the dimension of $H_*(\Mor_{\cA(T^2)}(\CFDh(S^3-U), \CFDh(S^3-J)))$, our next task is to find a basis. To simplify the exposition, we introduce some notation. The summand of $\cCFK^-(J)$ generated by $x$ gives rise to a summand of $\CFDh(S^3-J)$ generated by an element we will also denote $x$, with differential $\delta^1 (x) = \rho_{12}x$. Call this type-D structure $B$. The unit boxes in $\cCFK^-(J)$ correspond to boxes in $\CFDh(S^3-J)$ as in Figure \ref{fig: unit box type d structure}. Let $C$ be such a type-D structure.

\begin{figure}[h!]
        \begin{tikzcd}
            b \ar[d, "\rho_1"]
                & y^1 \ar[l, "\rho_2"]
                    & a \ar[l, "\rho_3"] \ar[d, "\rho_1"] \\
            y^2 
                &     
                    & y^4 
                         \\
            e \ar[u, "\rho_{123}"]
                & y^3 \ar[l, "\rho_2"]
                    & c \ar[l, "\rho_3"] \ar[u, "\rho_{123}"].
        \end{tikzcd}
        \caption{The type-D structure $C$ associated to a unit box.}\label{fig: unit box type d structure}
    \end{figure}

    If we ignore the gradings, $\CFDh(S^3-U)$ is isomorphic to $B$, and $\CFDh(S^3-J)$ is isomorphic to $B \oplus C^{\oplus 4}$. Therefore, it suffices to find bases for the homologies of $\Mor_{\cA(T^2)}(B, B)$ and $\Mor_{\cA(T^2)}(B, C)$.\\

    The space of homotopy classes of maps $B \ra B$ has a basis given by:
        \[
        \phi = (x \mapsto x)
        \]
        \[
        \psi = (x\mapsto \rho_{12}x).
        \]
    Since $H_*(\Mor_{\cA(T^2)}(B, B \oplus C^{\oplus 4})$ is ten dimensional, $H_*(\Mor_{\cA(T^2)}(B, C))$ must be two dimensional. By inspection, a basis is given by:
        \[
        g = (x \mapsto e + \rho_3 y^2 + \rho_1 y^3)
        \]
        \[
        h = (x \mapsto \rho_1 y^4).
        \]
    This is confirmed by Zhan's bordered Floer homology calculator \cite{bfh_python}.

\subsection{From concordance maps to complement maps}\label{subsection: concordance maps to complement maps}

The maps $F_{S^3\times I - C}$ and $F_{S^3\times I - C'}$ satisfy the property that 
\[
\bI_{\CFA^-(T_\infty, \lambda)} \boxtimes F_{S^3\times I - C} \sim F_C
\]
and 
\[
\bI_{\CFA^-(T_\infty, \lambda)} \boxtimes F_{S^3\times I - C'} \sim F_{C'}.
\]
Therefore, to determine these maps, we will compute $\bI_{\CFA^-(T_\infty, \lambda)} \boxtimes f$ for all basis elements, $f$, of the morphism space.

\begin{figure}[h!]
    \includegraphics[scale = .5]{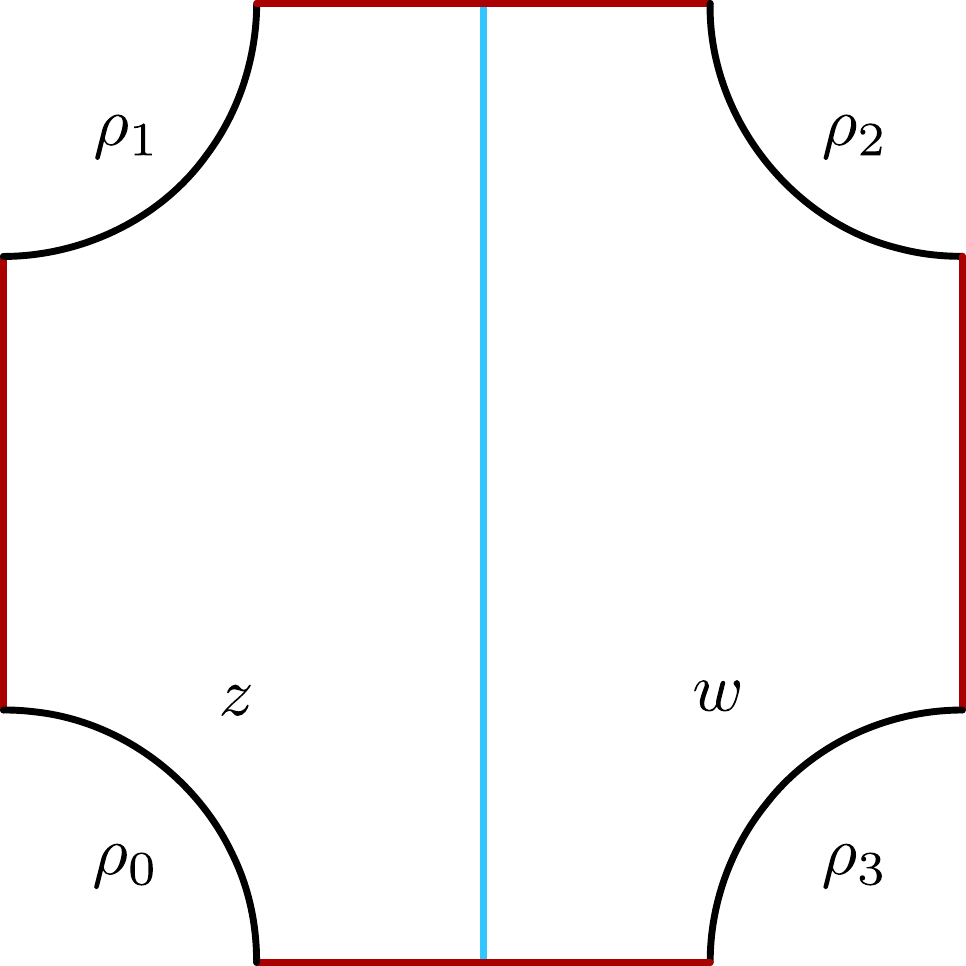}
    \caption{A doubly pointed bordered Heegaard diagram for the longitudinal unknot in the solid torus, $(T_\infty, \lambda)$.}\label{fig: doubly pointed heegaard diagram for unknot in solid torus.}
\end{figure}

$\CFA^-(T_\infty, \lambda)$ is a right $\cA_\infty$-module over $\cA(T^2)$, whose $\cA_\infty$-operations can be computed by counting holomorphic disks in the doubly pointed bordered Heegaard diagram shown in Figure \ref{fig: doubly pointed heegaard diagram for unknot in solid torus.}.

Let $\alpha$ be the single intersection point in $\alpha \cap \beta$. The only non-trivial $\cA_\infty$-operations are given by 
\[
m_{3 + j}(\alpha, \rho_3, \overbrace{\rho_{23}, \hdots, \rho_{23}}^{j}, \rho_2) = U^{j+1} \alpha.
\]

Call this module $A$. Since $A \boxtimes B$ has a single generator, $\alpha \otimes x$, the maps $\bI_{A} \boxtimes f$ are determined by the image of this element.

\begin{lemma}\label{lemma: tensor product for unknot in solid torus}
    Let $\phi, \psi, g, h$ be the basis for $H_*(\Mor_{\cA(T^2)}(B, B \oplus C))$ described above. Then, 
    \begin{align*}
        \bI_A \boxtimes \phi &= (\alpha \otimes x \mapsto \alpha \otimes x)\\
        \bI_A \boxtimes \psi &= 0\\
        \bI_A \boxtimes g &= (\alpha \otimes x \mapsto \alpha \otimes e)\\
        \bI_A \boxtimes h &= 0.
    \end{align*}
\end{lemma}

\begin{proof}
    The differentials of $C$ are shown again, as they are needed below.
    \begin{center}
        \begin{tikzcd}
            b \ar[d, "\rho_1"]
                & y^1 \ar[l, "\rho_2"]
                    & a \ar[l, "\rho_3"] \ar[d, "\rho_1"] \\
            y^2 
                &     
                    & y^4 
                         \\
            e \ar[u, "\rho_{123}"]
                & y^3 \ar[l, "\rho_2"]
                    & c \ar[l, "\rho_3"] \ar[u, "\rho_{123}"].
        \end{tikzcd}
    \end{center}
    Since many of the differentials on the tensor product are trivial, many terms will be forced to zero. In particular, $\delta^1 (x) = \rho_{12} x$, and since there are no $\cA_\infty$-operations involving $\rho_{12}$, any terms involving $\delta^1 (x)$ will be zero. $\bI_A \boxtimes \psi$ is zero for this reason. We make use of the graphical notion of \cite[Chapter 2]{LOT_bordered_HF}. 

    By strict unitality, $\bI_A \boxtimes \vphi(\alpha \otimes x)$ has a single term, $\alpha \otimes x$.

    \[\bI_{A} \boxtimes \phi = 
    \begin{tikzcd}[column sep = tiny]
        \alpha\ar[dd] & & x\ar[d] \\
          & & \phi \ar[dd, "x"] \ar[dll, "1"] \\
        m \ar[d, "\alpha"] & &  \\
        \alpha & & x
    \end{tikzcd} = (\alpha \otimes x \mapsto \alpha \otimes x), \]

    To compute $\bI_A \boxtimes g = \bI_A \boxtimes (x \mapsto e + \rho_3y^2 + \rho_1 y^3)$, we first note that, again, by strict unitality, the only term $e$ could contribute is $\alpha \otimes e$. Secondly, any term contributed by $\rho_1 y^3$ will be of the form $m_\ell(\alpha, \overbrace{\rho_{12}, \hdots, \rho_{12}}^i, \rho_1, \hdots )\otimes \xi$. But, since none of the $\cA_\infty$-operations involve $\rho_{1}$, any term of this form must be zero. Therefore, all that remains is to check whether $\rho_3 y^2$ contributes any nonzero terms to $\bI_A \boxtimes g$. A nonzero term could appear, since all the $\cA_\infty$-operations involve $\rho_3$, but $\delta^1(y^2) = 0$, so no $\rho_2$ coefficient will appear. Therefore, 
    \[
    \bI_{A} \boxtimes g = 
    \begin{tikzcd}[column sep = tiny]
        \alpha\ar[dd] & & x\ar[d] \\
          & & g\ar[dd, "e"] \ar[dll, "1"] \\
        m \ar[d, "\alpha"] & &  \\
        \alpha & & e
    \end{tikzcd} + 
    \begin{tikzcd}[column sep = tiny]
        \alpha\ar[dd] & & x\ar[d] \\
          & & g\ar[dd, "y^2"] \ar[dll, "\rho_3"] \\
        m \ar[d] & &  \\
        0 & & y^2
    \end{tikzcd} + 
    \begin{tikzcd}[column sep = tiny]
        \alpha\ar[ddd] & & x\ar[d] \\
          & & g\ar[d, "y^3"] \ar[ddll, "\rho_1" left] \\
          & & \vdots \ar[d] \ar[dll] \\
        m \ar[d] & & \delta^1 \ar[d] \ar[ll]  \\
        0 & & \xi
    \end{tikzcd} = (\alpha \otimes x \mapsto \alpha \otimes e)
    \]

    Consider $\bI_A \boxtimes h = \bI_A \boxtimes (x\mapsto \rho_1 y^4)$. Much like the previous case, $\rho_1 y^4$ can contribute no nonzero terms, since $\rho_1$ does not appear in any of the $\cA_\infty$-operations.
    \[\bI_{A} \boxtimes h = 
    \begin{tikzcd}[column sep = tiny]
        \alpha\ar[dd] & & v\ar[d] \\{}
          & & h\ar[dd, "y^4"] \ar[dll, "\rho_1"] \\
        m \ar[d, "0"] & &  \\
        \alpha & & y^4
    \end{tikzcd} = (\alpha \otimes v \mapsto 0). \]
    This concludes this computation.
\end{proof}

We can now turn to the maps $F_{S^3\times I - C}$ and $F_{S^3\times I - C'}$. Recall from Section \ref{Section: computations of exotic disk maps}:
\[
F_C(v)+ F_{C'}(v) = e_1 + e_2,
\]
where $v$ is the generator for $\HFK^-(U)$. A basis for $H_*(\Mor_{\cA(T^2)}(\CFDh(S^3-U), \CFDh(S^3-J)))$ is given by $\{\phi, \psi, g_i, h_i\}$, where $g_i$ and $h_i$ are maps from $\CFDh(S^3-U)$ to the $i$th box in $\CFDh(S^3-J))$ which agree with the maps $g$ and $h$ above. The complex $\CFDh(S^3-J)$ is shown in full in Figure \ref{figure: Full CFD}.

\begin{figure}[h!]
\begin{minipage}{.4\textwidth}
\[
    \begin{tikzcd}
            b_i \ar[d, "\rho_1"]
                & y^1_i \ar[l, "\rho_2"]
                    & a \ar[l, "\rho_3"] \ar[d, "\rho_1"] \\
            y^2_i 
                &     
                    & y^4_i 
                         \\
            e_i \ar[u, "\rho_{123}"]
                & y^3_i \ar[l, "\rho_2"]
                    & c_i \ar[l, "\rho_3"] \ar[u, "\rho_{123}"]
        \end{tikzcd}
\]
\end{minipage}
\begin{minipage}{.4\textwidth}
\[
\begin{tikzcd} g_i \ar[d, "\rho_1"]
            & z^1_i \ar[l, "\rho_2"]
                & f_i \ar[l, "\rho_3"] \ar[d, "\rho_1"] \\
        z^2_i 
            &     
                & z^4_i 
                    \\
        j_i \ar[u, "\rho_{123}"]
            & z^3_i \ar[l, "\rho_2"]
                & h_i \ar[l, "\rho_3"] \ar[u, "\rho_{123}"]
\end{tikzcd}
\begin{minipage}{.1\textwidth}
\begin{tikzcd}
    x  \ar[loop right, "\rho_{12}"] 
\end{tikzcd}
\end{minipage}
\]
\end{minipage}
\caption{The full complex $\CFDh(S^3-J)$.}\label{figure: Full CFD}
\end{figure}
Let $v$ also denote the single generator for $\CFDh(S^3-U)$. By Lemma \ref{lemma: tensor product for unknot in solid torus}, we can identify which maps $f: \CFDh(S^3-U) \ra \CFDh(S^3-J)$ have the property that $\bI_A \boxtimes f$ is homotopic to either $F_C$ or $F_{C'}$. \\

Again, $F_C(v) + F_{C'}(v) = e_1 + e_2$. Since the maps $\psi, h_1, \hdots, h_4$ satisfy $\bI_A \boxtimes \psi = \bI_A \boxtimes h_i = 0$, we cannot rule out the possibility that they appear as terms in $F_{S^3\times I - C}(v)$ and $F_{S^3\times I - C'}(v)$. Therefore, we can deduce that
\[
F_{S^3\times I - C} + F_{S^3\times I - C'} \sim g_1 + g_2 + \vep_1\cdot \psi + \vep_2 \cdot\sum_{i=1}^4 h_i
\]
where $\vep_i \in \{0, 1\}$. Surprisingly, despite the indeterminacy of these maps, this will be sufficient information to prove Theorem \ref{Theorem 1: exotic disks with large stabilization distance}. 

In summary, we have the following:

\begin{prop}\label{prop: candidates for the complement maps}
    Let $C$ and $C'$ be the exotic concordances from the unknot to $J$. Then, the maps $F_{S^3\times I - C}$ and $F_{S^3\times I - C'}$ satisfy:
    \[
    F_{S^3\times I - C} + F_{S^3\times I - C'} \sim g_1 + g_2 + \vep_1\cdot \psi + \vep_2 \cdot\sum_{i=1}^4 h_i
    \]
    where $\vep_i \in \{0, 1\}$.
\end{prop}

\begin{remark}
    It seems likely that one could show that $\vep_2$ and $\vep_4$ are actually zero by considering the grading on $\CFDh(S^3-J)$. However, it will turn out to be of no consequence, so we do not pursue this.
\end{remark}

\subsection{From complement maps to cabled concordance maps}\label{subsection: from complement maps to cabled concordance}

Having found our candidates for the maps $F_{S^3\times I - C}$ and $F_{S^3\times I - C'}$, all that remains is to compute the tensor products of the candidates with the identity map for the $\cA_\infty$-module associated to the $(p, 1)$-cable pattern in the solid torus.

A doubly pointed Heegaard diagram $\cH_P$ for the $(p, 1)$-cable in the solid torus is shown in Figure \ref{fig: (p, 1) cable diagram}.
\begin{figure}[h!]
    \centering
    \includegraphics[scale=.5]{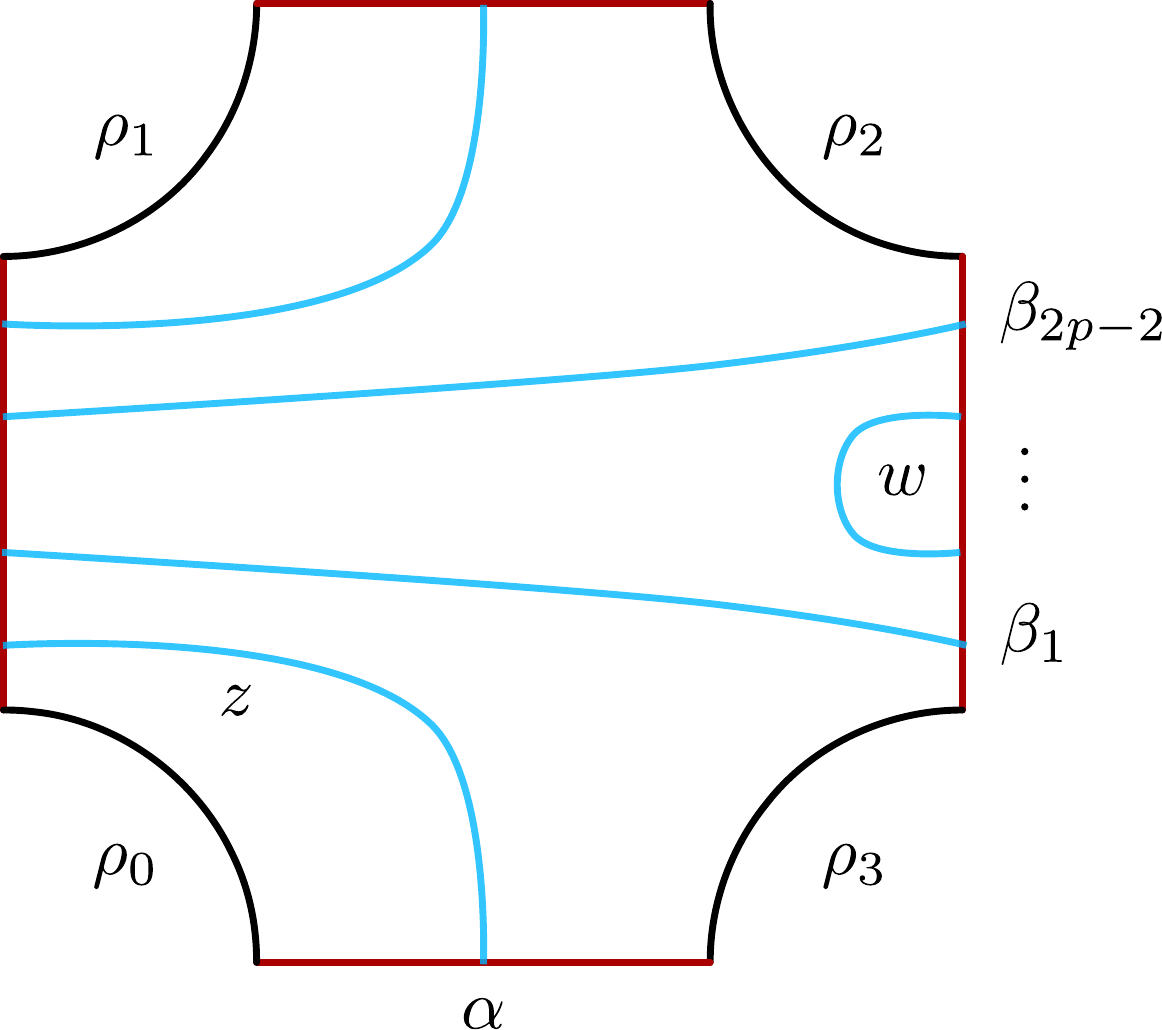}
    \caption{A doubly pointed bordered Heegaard diagram for the $(p, 1)$-cable in the solid torus.}
    \label{fig: (p, 1) cable diagram}
\end{figure}
Let $A_p = \CFA^-(\cH_p)$. $A_p$ is generated by $\alpha, \beta_1, \hdots, \beta_{2p-2}$. Since we are only interested in computing maps $\CFA^-(\cH_p) \boxtimes \CFDh(S^3-U) \ra \CFA^-(\cH_p) \boxtimes \CFDh(S^3-J)$ and the homology of $\CFA^-(\cH_p)\boxtimes \CFDh(S^3-U) \simeq \CFK^-(U)$ is generated by $\alpha \otimes v$, it is enough to consider the $\cA_\infty$-operations involving $\alpha$. 

\begin{align*}
   m_{2+i}(\alpha, \overbrace{\rho_{12} , \hdots, \rho_{12}}^{i}, \rho_1) = \beta_{2p-i-2} &\hspace{1cm} 0 \le i \le p-2 \\
   m_{4+i+j}(\alpha, \rho_3, \overbrace{\rho_{23}, \hdots, \rho_{23}}^{j}, \rho_2, \overbrace{\rho_{12}, \hdots, \rho_{12}}^{i}, \rho_1) = U^{pj+i+1}\beta_{i+1} &\hspace{1cm} 0 \le i \le p-2, 0 \le j\\
   m_{3+i}(\alpha, \rho_3, \overbrace{\rho_{23}, ..., \rho_{23}}^{j}, \rho_2) = U^{p(j+1)}\alpha, &\hspace{1cm} 0 \le j.
\end{align*}

For the full collection of $\cA_\infty$-operations, see \cite[Section 4]{petkova_cables_thin_knots}. As before, we will start by computing $\bI_{A_p} \boxtimes f$ for the basis elements of $H_*(\Mor_{\cA(T^2)}(B, B \oplus C))$, and then use the fact that $\CFDh(S^3-J)$ is isomorphic to $B \oplus C^{\oplus 4}$ to compute $F_{S^3\times I - C_p}$ and $F_{S^3\times I - C_p'}$.

\begin{lemma}\label{lemma: cable box tensor product}
    Let $\phi$, $\psi$, $h$, and $h$ be the basis of $H_*(\Mor_{\cA(T^2)}(B, B \oplus C))$ as computed in Section \ref{subsection: concordance maps to complement maps}. Then, 
    \begin{align*}
        \bI_{A_p} \boxtimes \phi &= (\alpha \otimes x \mapsto \alpha \otimes x) \\
        \bI_{A_p} \boxtimes \psi &= (\alpha \otimes x \mapsto 0) \\
        \bI_{A_p} \boxtimes g &= \left(\alpha \otimes x \mapsto \alpha \otimes e + \sum_{i=0}^{p-2} \beta_{2p-i-2}\otimes y^3\right) \\
        \bI_{A_p} \boxtimes h &= \left(\alpha \otimes x \mapsto  \sum_{i=0}^{p-2} \beta_{2p-i-2}\otimes y^4\right)
    \end{align*}
\end{lemma}
\begin{proof}
    We proceed as in Lemma \ref{lemma: tensor product for unknot in solid torus}. This computation is slightly more involved, as there are more $\cA_\infty$-operations. 

    First, for the map $\psi = (x \mapsto \phi_{12}x)$, it must be that $\bI_{A_p} \boxtimes \psi = 0$. Since $\delta^1(x) = \rho_{12} x$, any term coming from $\rho_{12}x$ will be of the form $m_{k}(\rho_{12}, \hdots, \rho_{12})$, which must be zero, as there are no $\cA_\infty$-operations only involving $\rho_{12}$. 

    For the map $\phi = (x \mapsto x)$, the map $\bI_{A_p} \boxtimes \phi$ has a single term by strict unitality:
    \[
    \bI_{A_p} \boxtimes \phi = 
    \begin{tikzcd}[column sep = tiny]
        \alpha\ar[dd] & & v\ar[d] \\
          & & f\ar[dd, "v"] \ar[dll, "1"] \\
        m \ar[d, "\alpha"] & &  \\
        \alpha & & v
    \end{tikzcd} = (\alpha \otimes v \mapsto \alpha \otimes v).
    \]
    The map $h = (x \mapsto \rho_1 y^4)$ has many nonzero terms, since there are nontrivial $\cA_\infty$-operations involving $\rho_{12}$ and $\rho_1$, namely $m_{2+i}(\alpha, \overbrace{\rho_{12} , \hdots, \rho_{12}}^{i}, \rho_1) = \beta_{2p-i-2}$ for $0 \le i \le p-2$. But these are the only terms that can appear, since $\delta^1(y^4) = 0$. Therefore, 
    \[
    \bI_{A_p} \boxtimes h = 
    \begin{tikzcd}[column sep = tiny]
        \alpha\ar[dddd] & & v\ar[d] \\
          & & \delta^1\ar[d, "v"] \ar[dddll, "\rho_{12}" left] \\
          & & \vdots \ar[d, "v"] \ar[ddll, "\rho_{12}" right]\\
          & & g \ar[dd, "y^4"] \ar[dll, "\rho_1"]\\
         m \ar[d, "\beta_{2p-i-2}"]& &  \\
         \beta_{2p-i-2} & & y^4
    \end{tikzcd} = \left(\alpha \otimes v \mapsto  \sum_{i=0}^{p-2} \beta_{2p-i-2}\otimes y^4\right).
    \]

    The map $g = (x \mapsto e + \rho_3 y^2 + \rho_1 y^3)$ is the most complicated. Once again, by strict unitality, it must be that $\alpha \otimes e$ appears as a term in $\bI_B \boxtimes g(\alpha \otimes x)$, and no other terms will come from $e$. The differential of $y^2$ is zero, so $\rho_2 y^2$ can only contribute terms of the form $m_{2+i}(\alpha, \overbrace{\rho_{12}, \hdots, \rho_{12}}^{i}, \rho_2) \otimes y^2$, but these are all zero. Finally, the differential of $y^3$ is $\rho_2 e$, and the differential of $e$ is $\rho_{123} y^2$, so there could potentially be terms of the form $m_{2+i}(\alpha, \overbrace{\rho_{12}, \hdots, \rho_{12}}^{i}, \rho_1) \otimes y^3$,  $m_{3+i}(\alpha, \overbrace{\rho_{12}, \hdots, \rho_{12}}^{i}, \rho_1, \rho_2)\otimes e$ or $m_{4+i}(\alpha, \overbrace{\rho_{12}, \hdots, \rho_{12}}^{i}, \rho_1, \rho_2, \rho_{123})\otimes y^2$, but only the first case introduces non-zero terms. Therefore, 
    \[
    \bI_{A_p} \boxtimes g = 
    \begin{tikzcd}[column sep = tiny]
        \alpha\ar[dd] & & v\ar[d] \\
          & & g\ar[dd, "e"] \ar[dll, "1"] \\
        m \ar[d, "\alpha"] & &  \\
        \alpha & & e 
    \end{tikzcd} +
    \begin{tikzcd}[column sep = tiny]
        \alpha\ar[dddd] & & v\ar[d] \\
          & & \delta^1\ar[d, "v"] \ar[dddll, "\rho_{12}" left] \\
          & & \vdots \ar[d, "v"] \ar[ddll, "\rho_{12}" right]\\
          & & g \ar[dd, "y^3"] \ar[dll, "\rho_1"]\\
         m \ar[d, "\beta_{2p-i-2}"]& &  \\
         \beta_{2p-i-2} & & y^3
    \end{tikzcd} 
    = \left(\alpha \otimes v \mapsto \alpha \otimes e + \sum_{i=0}^{p-2} \beta_{2p-i-2}\otimes y^3\right).
    \]
    This concludes the computation.
\end{proof}

In Proposition \ref{prop: candidates for the complement maps}, we found the sum of the maps $F_{S^3\times I - C}$ and $F_{S^3\times I - C'}$. By Theorem \ref{Theorem 2: satellite formula for concordance maps}, the maps $F_{C_p}$ and $F_{C_p'}$ induced by the $(p, 1)$-cables of $C$ and $C'$ can be computed as $\bI_{A_p} \boxtimes F_{S^3\times I - C}$ and $\bI_{A_p} \boxtimes F_{S^3\times I - C'}$ respectively. Lemma \ref{lemma: cable box tensor product} can now be applied to give us the candidates for the maps $F_{C_p}$ and $F_{C_p'}$. As before, let $\{\phi, \psi, g_i, h_i\}$ be the basis for $H_*(\Mor_{\cA(T^2)}(\CFDh(S^3-U), \CFDh(S^3-J)))$ from Section \ref{subsection: concordance maps to complement maps}. 

Since 
    \[
    F_{S^3\times I - C} + F_{S^3\times I - C'}\sim g_1 + g_2 + \vep_1\psi + \vep_2 \cdot\sum_{i=1}^4 h_i
    \]
for $\vep_i \in \{0, 1\}$, from Lemma \ref{lemma: cable box tensor product} it follows that
    \[
        (\bI_{A_p} \boxtimes F_{S^3\times I - C} + \bI_{A_p} \boxtimes F_{S^3\times I - C'})(\alpha \otimes v) = 
    \]
    \[   \alpha\otimes (e_1 +e_2 ) + \left( \sum_{k=1}^2 \sum_{i=0}^{p-2} \beta_{2p-i-2} \otimes y^3_k + \vep_2 \cdot \beta_{2p-i-2} \otimes (y^4_k + z^4_k) \right) 
    \]
    in homology. Summarizing our results, we have the following.

\begin{prop}\label{prop: candidates for C_p and C_p'}
    Let $C_p$ and $C_p'$ be the $(p, 1)$-cables of the exotic concordances $C$ and $C'$. Let $\cH_p$ be the doubly pointed bordered Heegaard diagram for the $(p, 1)$-cable pattern shown in Figure \ref{fig: (p, 1) cable diagram}. Then, the maps $F_{C_p}$ and $F_{C_p'}$ satisfy:
    \[
        (\bI_{A_p} \boxtimes F_{S^3\times I - C} + \bI_{A_p} \boxtimes F_{S^3\times I - C'})(\alpha \otimes v) = 
    \]
    \[   \alpha\otimes (e_1 +e_2 ) + \left( \sum_{k=1}^2 \sum_{i=0}^{p-2} \beta_{2p-i-2} \otimes y^3_k + \vep_2 \cdot \beta_{2p-i-2} \otimes (y^4_k + z^4_k) \right) 
    \]
    where $\vep$ is either 0 or 1.
\end{prop}


\section{Disks with large stabilization distance}\label{Section: large stabilization distance}

\subsection{Stabilization distance bounds}

We are now equipped to prove Theorem \ref{Theorem 1: exotic disks with large stabilization distance}. $\CFK^-(J_{(p, 1)})$ can be computed by, again, appealing to the work of Petkova \cite{petkova_cables_thin_knots}, who computes the differential on complexes of the form $\CFA^-(\cH_p) \boxtimes C$ where $C$ is a unit box as in Figure \ref{fig: unit box type d structure}. The subcomplex of $\CFK^-(J_{(p, 1)})$ corresponding to the unit boxes generated by $\{a_i, b_i, c_i, e_i\}$ is shown. 

\begin{center}
    \begin{tikzcd}[column sep = small]
        \alpha \otimes c_k \ar[d, "U^p"] & 
            \beta_i \otimes y^3_k \ar[d, "U^{p-i}"] & 
                \beta_j \otimes y^4_k \ar[d, "U^{p-j}"] & 
                    \beta_{p-1} \otimes y_k^1 \ar[d, "U"] & 
                        \beta_{p-1} \otimes y_k^3 \ar[d, "U"]\\
        \alpha\otimes e_k & 
            \beta_{2p-i-1} \otimes y_k^3 &
                \beta_{2p-j-1} \otimes y_k^4 & 
                    \beta_p \otimes y_k^1 & 
                        \beta_p \otimes y_k^3 \\
        \beta_1 \otimes y_k^4 \ar[d, "U^{p-1}"]& 
            \alpha \otimes a_k \ar[r, "U^p"], \ar[d, "U"] \ar[dl] &
                \alpha \otimes b_k \ar[d] & 
                    \beta_i \otimes y_k^1 \ar[d, "U^{p-i}"] \ar[r, "U"] &
                        \beta_{i+1} \otimes y_k^2 \ar[d,"U^p-i-1"]\\ 
        \beta_{2p-2} \otimes y^4_k  & 
            \beta_1 \otimes y_k^2 \ar[r, "U^{p-1}"] &
                \beta_{2p-2} \otimes y_k^2 & 
                    \beta_{2p-i-1} \otimes y_k^1 \ar[r] & 
                        \beta_{2p-i-1}\otimes y_k^2 
    \end{tikzcd}
\end{center}
Here, $k = 1, 2$, $1 \le i \le p-2$, and $2 \le j \le p-2$. The subcomplex of $\CFK^-(J_{(p, 1)})$ corresponding to the unit boxes generated by $\{f_i, g_i, h_i, j_i\}$ is obtained by replacing $y_k^i$ by $z_k^i$ and $a_i, b_i, c_i, e_i$ with $f_i, g_i, h_i, j_i.$ 

\begin{prop}\label{prop: tau for D_p D_p'}
    For the $(p, 1)$-cables of the exotic disks $D$ and $D'$,
    \[
        \tau(D_p, D_p') = p.
    \]
    Therefore, the stabilization distance between $D_p$ and $D_p'$, $\mu(D_p, D_p')$ is at least $p$.
\end{prop}
\begin{proof}
By inspecting the complex $\CFK^-(J_{p, 1})$, we make the following observations: the $\F[U]$-tower in $\HFK^-(J_{p, 1})$ is generated by $\alpha \otimes v$; $\alpha \otimes e_1$ and $\alpha \otimes e_2$ both have torsion order $p$ and are distinct in $\HFK^-(J_{(p,1)})$; the terms $b_{2p-i-2} \otimes y_k^3$, $b_{2p-i-2} \otimes y_k^4$, and $b_{2p-i-2} \otimes z_k^4$ have torsion order at most $(p-1)$.

By Proposition \ref{prop: candidates for C_p and C_p'}, 
    \[
        (F_{C_p} + F_{C_p'})(\alpha \otimes v) = \alpha\otimes (e_1 +e_2 ) + \left( \sum_{k=1}^2 \sum_{i=0}^{p-2} \beta_{2p-i-2} \otimes y^3_k + \vep_2 \cdot \beta_{2p-i-2} \otimes (y^4_k + z^4_k) \right) 
    \]
By considering the torsion orders of the elements in these sums, we see that as elements of $\HFK^-(J_{(p, 1)})$,
\[
U^{p-1}\cdot (F_{C_p}( \alpha \otimes v) + F_{C_p'}( \alpha \otimes v)) = U^{p-1}\cdot\alpha \otimes (e_1 + e_2) \ne 0,
\]
which implies $\tau(D_p, D_p') \ge p$. Since $\alpha \otimes (e_1 + e_2)$ has torsion order $p$, we in fact have that  $\tau(D_p, D_p') = p$. Therefore, by the Juh\'asz-Zemke lower bound for the stabilization distance, we have that 
\[
    p = \tau(D_p, D_p') \le \mu(D_p, D_p').
\]
\end{proof}

Combining this result with Proposition \ref{prop: upper bound for cabled stabilization distance} proves Theorem \ref{Theorem 1: exotic disks with large stabilization distance}.

\section{Exotic cobordisms}

Let $C: K_1 \ra K_2$ be a concordance. Following \cite{gordon_knots_homology_spheres}, the concordance $C$ can be used to construct a cobordism between $S_n(K_1)$ and $S_n(K_2)$ as follows. Remove a tubular neighborhood of $C$ from $S^3 \times I$. The framing of $K_1$ determines a map $h: S^1\times \partial D^2 \ra \partial \nu(K_1)$. Define the cobordism $W_n$ as $(S^1\times D^2 \times I) \cup_{h \times \id} (S^3\times I-C)$. Since $W_n$ is built by gluing $S^1\times D^2 \times I$ to a concordance complement, the associated map $\wh{F}_{W_n}$ can be computed using Theorem \ref{Theorem 2: satellite formula for concordance maps}, as a tensor product $\bI_{\CFAh(S^1\times D^2, n)} \boxtimes F_{S^3\times I-C}$, where $\CFAh(S^1\times D^2, n)$ is the type-A structure associated to the $n$-framed solid torus. \\

In particular, by taking the exotic concordances $C, C'$ from the unknot to $J$, this construction yields a pair of cobordisms $W_n, W_n': L(n, 1) \ra S_n^3(J)$. The ambient topological isotopy from $C$ to $C'$ together with the existence of topological normal bundles implies $W_n$ and $W_n'$ are homeomorphic. However, as we will show, $W_n$ and $W_n'$ are not diffeomorphic rel boundary, as they induce different maps on $\HFh$.

\begin{figure}
    \includegraphics[scale=.5]{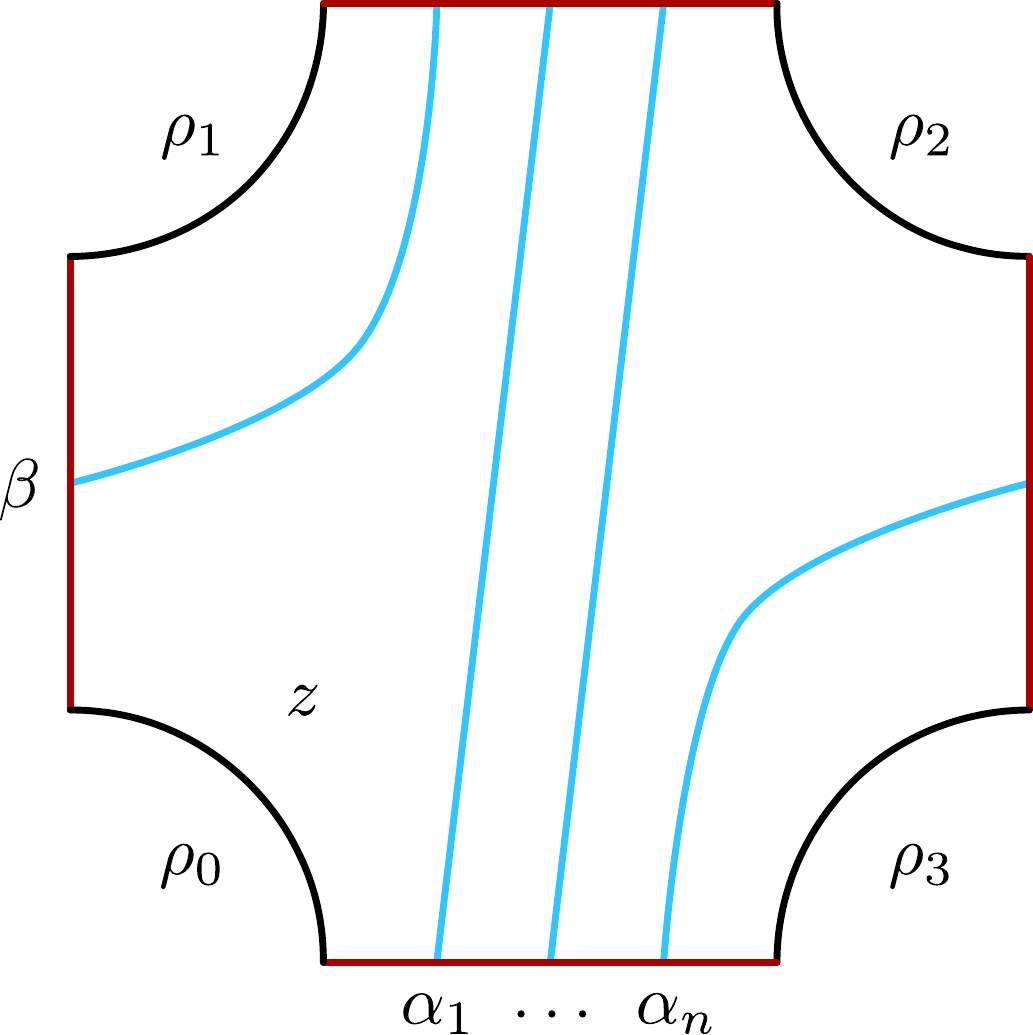}
    \caption{A bordered Heegaard diagram for the $n$-framed solid torus.}
    \label{fig: n framed torus}
\end{figure}
\begin{proof}[Proof of Theorem \ref{theorem: exotic cobordisms}]
    
A bordered Heegaard diagram for the $n$-framed solid torus is shown in Figure \ref{fig: n framed torus}.  The nontrivial $\cA_\infty$-operations on $\CFAh(\cH_n)$ are given by 
\begin{align*}
m_2(\alpha_1, \rho_1) = \beta, & \\
m_3(\alpha_i, \rho_3, \rho_2)  = \alpha_{i+1}, & \hspace{1cm} 1 \le i \le n-1, \\
m_2(\alpha_n, \rho_3) = \beta.  &   
\end{align*}

The knot $J$ has genus 1, so for $n \ge 1$, $\HFh(S_n^3(J), [s])$ can be computed using the large surgery formula:
\[
\HFh(S^3_n(J), [s]) \cong \begin{cases}
    \F^{9} & s = 0\\ 
    \F & 1 \le s \le n-1.
\end{cases}
\]

However, as the maps $\wh{F}_{W_n}$ and $\wh{F}_{W_n'}$ will be computed as a box tensor product, it is convenient to compute $\HFh(S^3_n(J))$ using the pairing theorem. Recall that $\CFDh(S^3-J)$ is shown in Figure \ref{figure: Full CFD}. The portion of $\CFAh(\cH_n) \boxtimes_{\cA(T^2)} \CFDh(S^3-J)$ corresponding to the unit boxes in $\cCFK^-(J)$ generated by $\{a_i, b_i, c_i, e_i\}$ and the generator $x$ is shown below.

\begin{center}
    \begin{tikzcd}
        & \alpha_1 \otimes c_i \ar[d] & \hdots & \alpha_{n-1} \otimes c_i \ar[d]& \alpha_n \otimes c_i \ar[d] \\
    \alpha_1 \otimes e_i & \alpha_2 \otimes e_i & \hdots & \alpha_n \otimes e_i & \beta \otimes y^3_i \\
    \end{tikzcd}
    \begin{tikzcd}
       \beta \otimes y^4_i & \alpha_1 \otimes a_i \ar[d]\ar[l] & \hdots & \alpha_{n-1} \otimes a_i \ar[d]& \alpha_n \otimes a_i \ar[d] \\
    \alpha_1 \otimes b_i \ar[d] & \alpha_2 \otimes b_i & \hdots & \alpha_n \otimes b_i & \beta \otimes y^1_i \\
    \beta \otimes y^2_i & \alpha_1 \otimes x & \hdots & \alpha_{n-1} \otimes x & \alpha_n \otimes x 
    \end{tikzcd}
\end{center}
The rest of the complex is obtained by replacing $a_i, b_i, c_i, e_i$ with $f_i, g_i, h_i, j_i$ and $y_i^k$ with $z_i^k$. The homology of this complex is generated by $\alpha_1\otimes e_i$, $\alpha_2 \otimes b_i = \beta \otimes y_i^4$, $\alpha_1 \otimes j_i$, $\alpha_2 \otimes g_i = \beta \otimes z_i^4$ and $\alpha_1 \otimes x, \hdots, \alpha_n \otimes x$. So, this computation does agree with the large surgery computation. \\

The maps $\bI_{\CFAh(\cH_n)} \boxtimes F_{S^3\times I - C}$ and $\bI_{\CFAh(\cH_n)} \boxtimes F_{S^3\times I - C'}$ can now be computed. In Proposition \ref{prop: candidates for the complement maps}, we showed that
\[
F_{S^3\times I - C} + F_{S^3\times I - C'} \sim g_1 + g_2 + \vep_1\cdot\psi + \vep_2 \cdot\sum_{i=1}^4 h_i.
\]

Recall that $\phi = (v \mapsto x)$, $\psi = (v \mapsto \rho_{12} x)$, $g_i = (v \mapsto e_i + \rho_3 y^2_i + \rho_1 y^3_i)$, and $h_i = (v \mapsto \rho_1 y^4_i)$.

First, we compute $\bI_{\CFAh(\cH_n)} \boxtimes \phi$. By strict unitality, we have that 
\[
\bI_{\CFAh(\cH_n)} \boxtimes \phi =
    (\alpha_k \otimes v \mapsto \alpha_k \otimes x), \;\; 1 \le k \le n.
\]
Since $\delta^1 (v) = \rho_{12}v$ and $\delta^1 (x) = \rho_{12} x$, and there are no non-trivial $\cA_\infty$-operations involving $\rho_{12}$, it must be that $\bI_{\CFAh(\cH_n)} \boxtimes \psi = 0.$

Now for $\bI_{\CFAh(\cH_n)} \boxtimes g_i$. The only term in $\bI_{\CFAh(\cH_n)} \boxtimes g_i(\alpha_k \otimes v)$ involving $e_i$ is $\alpha_k \otimes e_i$. Since $m_2(\alpha_n, \rho_3) = \beta$, we have that $\rho_3 y_i^2$ contributes a term to $\bI_{\CFAh(\cH_n)} \boxtimes g_i(\alpha_n \otimes v)$. The coefficient of $y_i^2$ of $\rho_3$, so there might have been a term coming from $m_3(\alpha_k, \rho_3, \rho_2)$, but since $\delta^1(y_i^2) = 0$, this is not the case.

Similarly, $m_2(\alpha_1, \rho_1) = \beta$, which means $\rho_1y^3_i$ contributes a term to $\bI_{\CFAh(\cH_n)} \boxtimes g_i(\alpha_1 \otimes v)$. Summarizing,
\[
\bI_{\CFAh(\cH_n)} \boxtimes g_i = \begin{cases}
    (\alpha_k \otimes v \mapsto \alpha_k \otimes e_i + \beta \otimes y_i^3) & k = 1\\
    (\alpha_k \otimes v \mapsto \alpha_k \otimes e_i) & 1 < k < n \\
    (\alpha_k \otimes v \mapsto \alpha_k \otimes e_i + \beta \otimes y_i^2) & k= n 
\end{cases}
\]

Lastly, we compute $\bI_{\CFAh(\cH_n)} \boxtimes h_i$. $\delta^1 (y^4_i) = 0$, so contributes nothing. There is an $\cA_\infty$-operation involving $\rho_1$, namely $m_2(\alpha_1, \rho_1) = \beta$, so the only non-trivial term is $\alpha_1 \otimes v \mapsto \beta \otimes y_i^4$. Altogether,

\[
\bI_{\CFAh(\cH_n)} \boxtimes h_i = \begin{cases}
    (\alpha_k \otimes v \mapsto \beta \otimes y^4_i) & k = 1 \\
    (\alpha_k \otimes v \mapsto 0) & k > 1 
\end{cases}
\]

Note that $\alpha_2 \otimes e_i, \hdots \alpha_n\otimes e_i$, $\beta \otimes y_i^2$, and $\beta \otimes y_i^3$ are nullhomologous, so as maps $\HFh(L(n, 1)) \ra \HFh(S_n^3(J))$, we have that 
\[
\wh{F}_{W_n}(\alpha_k \otimes v) + \wh{F}_{W_n'}(\alpha_k \otimes v) = \begin{cases}
    \alpha_k \otimes (e_1 + e_2) + \vep \cdot \beta \otimes (y_1^4 + y_2^4 + z_1^4 + z_2^4) & k = 1\\
   0 & \text{else}.
\end{cases}
\]
Therefore, the maps on homology are distinct, implying $W_n$ and $W_n'$ are exotic.
\end{proof}

\bibliographystyle{alpha} 
\bibliography{mathbib}
\vspace{1cm}

\end{document}